\documentclass[a4paper]{amsart}
\usepackage[hypertex]{hyperref}
\usepackage{amsmath,amssymb}
\usepackage{bbold}
\usepackage{psfrag}
\usepackage{graphicx}
\usepackage{subfigure}
\usepackage{leftidx}
\usepackage[all]{xy}
\usepackage[latin1]{inputenc}

\newtheorem{thm}{Theorem}[section]
\newtheorem{cor}[thm]{Corollary}

\newtheorem{prop}[thm]{Proposition}

\theoremstyle{definition}

\newtheorem{rem}[thm]{Remark}
\newtheorem{exa}[thm]{Example}

\newcommand{\co}{\colon}
\newcommand{\id}{\mathrm{id}}
\newcommand{\diag}{\mathcal{D}iag}

\newcommand{\rstl}{\mathrm{RSL}}

\newcommand{\opp}{\mathrm{op}}

\newcommand{\cc}{\mathcal{C}}
\newcommand{\bb}{\mathcal{B}}
\newcommand{\dd}{\mathcal{D}}

\newcommand{\uu}{\mathcal{U}}

\newcommand{\zz}{\mathcal{Z}}

\newcommand{\kk}{\Bbbk}
\newcommand{\kt}{$\Bbbk$\nobreakdash-\hspace{0pt}}
\newcommand{\trait}{\nobreakdash-\hspace{0pt}}
\newcommand{\Rt}{$\mathrm{R}$\nobreakdash-\hspace{0pt}}
\newcommand{\ti}{\mbox{-}\,}

\newcommand{\un}{\mathbb{1}}

\newcommand{\Ob}{\mathrm{Ob}}

\newcommand{\End}{\mathrm{End}}
\newcommand{\Hom}{\mathrm{Hom}}
\newcommand{\Nat}{{\textsc{Hom}}}

\newcommand{\vect}{\mathrm{vect}}
\newcommand{\rep}{\mathrm{rep}}

\newcommand{\tr}{\mathrm{tr}}

\newcommand{\rank}{\mathrm{rank}}

\newcommand{\lev}{\mathrm{ev}}
\newcommand{\rev}{\widetilde{\mathrm{ev}}}
\newcommand{\lcoev}{\mathrm{coev}}
\newcommand{\rcoev}{\widetilde{\mathrm{coev}}}
\newcommand{\ldual}[1]{\leftidx{^\vee}{\!#1}{}}
\newcommand{\rdual}[1]{{#1}^\vee}

\newcommand{\lldual}[1]{\leftidx{^{\vee\vee}}{\!#1}{}}

\newcommand{\rrdual}[1]{{#1}^{\vee\vee}}

\newcommand{\scaledraw}[1]{\includegraphics[scale=.8]{#1.eps}}
\newcommand{\scaleraisedraw}[2]{\raisebox{-#1\height}{\includegraphics[scale=.8]{#2.eps}}}
\newcommand{\sdraw}[3]{\raisebox{-#1\height}{\scalebox{#2}{\includegraphics{#3.eps}}}}
\providecommand{\bysame}{\leavevmode\hbox to3em{\hrulefill}\thinspace}

\begin{document}

\title{Categorical centers and Reshetikhin-Turaev invariants}
\author[A. Brugui\`eres]{Alain Brugui\`eres}
\author[A. Virelizier]{Alexis Virelizier}
\thanks{The second author thanks the organizers of the \emph{International Conference
on Quantum Topology} held at the Institute of Mathematics of Hanoi (Vietnam) in August 2007. This paper is an enhanced
version of the talk he gave there.} \email{bruguier@math.univ-montp2.fr \and virelizi@math.univ-montp2.fr}
\date{\today}

\begin{abstract}
A theorem of M\"uger asserts that the center $\zz(\cc)$ of a spherical \kt linear category $\cc$ is a modular category if
$\kk$ is an algebraically closed field and the dimension of $\cc$ is invertible. We generalize this result to the case where
$\kk$ is an arbitrary commutative ring, without restriction on the dimension of the category. Moreover we construct the
analogue of the Reshetikhin-Turaev invariant associated to $\zz(\cc)$ and give an algorithm for computing this invariant in
terms of certain explicit morphisms in the category $\cc$. Our approach is based on (a)~Lyubashenko's construction of the
Reshetikhin-Turaev invariant in terms of the coend of a ribbon category; (b)~an explicit algorithm for computing this
invariant via Hopf diagrams;~(c) an algebraic interpretation of the center of $\cc$ as the category of modules over a
certain Hopf monad $Z$ on the category $\cc$; (d)~a generalization of the classical notion of Drinfeld double to Hopf
monads, which, applied to the Hopf monad $Z$, provides an explicit description of the coend of $\zz(\cc)$ in terms of the
category $\cc$.
\end{abstract}
\maketitle

\setcounter{tocdepth}{1} \tableofcontents

\section*{Introduction}
In the early 90's, two new `quantum' invariants of $3$-manifolds were introduced: the \emph{Reshetikhin-Turaev invariant},
and the \emph{Turaev-Viro invariant}. The definition of the Reshetikhin-Turaev invariant $\mathrm{RT}_\bb$~\cite{RT2,Tur2}
involves a \emph{modular category}~$\bb$, that is, a ribbon fusion category over a commutative ring $\kk$ satisfying a
non-degeneracy condition (invertibility of the $S$-matrix). The algorithm for computing its value on a $3$-manifold consists
in presenting the manifold by surgery along a ribbon link and then taking a linear combination of colorings of this link by
simple objects of $\bb$.

Similarly, the definition of the Turaev-Viro invariant $\mathrm{TV}_\cc$ \cite{TV}, as revisited by Barrett and Westbury
\cite{BW}, involves a \emph{spherical category}, that is, a sovereign fusion category over a commutative ring $\kk$
such that left and right traces coincide.  The dimension $\dim\cc$ of $\cc$ (which is the
sum of squares of dimensions of simple objects) is moreover assumed to be invertible in~$\kk$. The algorithm for computing
$\mathrm{TV}_\cc(M)$ consists in presenting the $3$\trait manifold $M$ by a triangulation, coloring the edges of the
triangulation with simple objects of $\cc$, and then evaluating the colored tetrahedra  by means of the $6j$-symbols of
$\cc$.

If $\bb$ is a modular category, then it is also a spherical category, and the Reshetikhin-Turaev and Turaev-Viro invariants
are related \cite{Tur2,JR1} by:
\begin{equation*}
\mathrm{TV}_\bb(M)=\mathrm{RT}_\bb(M)\mathrm{RT}_\bb(-M)
\end{equation*}
for any $3$\ti manifold $M$, where $-M$ is the $3$-manifold $M$ with opposite orientation.

But in general a spherical category need not to be braided and so cannot be used as input to define the Reshetikhin-Turaev invariant.
However, spherical and modular categories are related by a theorem of M\"uger \cite{Mueg}: if~$\cc$ is a spherical
fusion category over an algebraically closed field $\kk$ and has invertible dimension, then its center $\zz(\cc)$ is a
modular fusion category of dimension  $\dim\zz(\cc)=(\dim\cc)^2$. In this setting, Turaev conjectured that, for
any $3$-manifold $M$,
\begin{equation*}
\mathrm{TV}_\cc(M)= \mathrm{RT}_{\zz(\cc)}(M).
\end{equation*}
This conjecture was shown to be true for some spherical categories~$\cc$ arising from subfactors, see \cite{KSW}. The
general case is still open.

In this context, a natural question is: how can we compute $\mathrm{RT}_{\zz(\cc)}(M)$? Using the algorithm given by
Reshetikhin and Turaev is not a practicable approach here, as that would require a description of the simple objects of
$\zz(\cc)$ in terms of those of $\cc$, and no such description is available in general. What we need is a different
algorithm for computing $\mathrm{RT}_{\zz(\cc)}(M)$, which one should be able to perform inside~$\cc$, without reference to
the simple objects of $\zz(\cc)$. This is the primary objective of this paper.

\medskip
In order to fulfill this objective, it will be convenient to adopt an alternative approach for constructing
$\mathrm{RT}$-like quantum invariants of $3$\trait manifolds, due to Lyubashenko \cite{Lyu2} and later developed in
\cite{LyuKer,Vir}, where the input data is a  (non-necessarily linear neither semisimple)  ribbon category $\bb$ which
admits a coend $C=\int^{X \in \bb} \ldual{X} \otimes X$. This coend $C$ is naturally endowed with a very rich algebraic
structure. In particular, it is a Hopf algebra in the braided category $\bb$ and comes equipped  with a Hopf pairing $\omega
\co C\otimes C \to \un$. Such a category $\bb$ is modular if the pairing $\omega$ is non-degenerate (this is the natural way
of formulating the invertibility of the $S$-matrix in this setting).

The construction of the Lyubashenko invariant consists in presenting the $3$\trait manifold by surgery along a ribbon link $L$,
using the universal property of the coend $C$ to associate a form $\phi_L$ to the link, and then evaluating this form on
an integral~$\Lambda$ of the Hopf algebra~$C$. Note that, more generally, one can evaluate the form $\phi_L$ by a `Kirby
element' $\alpha$ of $\bb$ to get other invariants $\tau_\bb(M;\alpha)$ of $3$-manifold invariants, see \cite{Vir}.
In particular, up to normalization, $\tau_\bb(M;\Lambda)$ is the Lyubashenko invariant and, in the special case where $\bb$ is a modular fusion category, $\tau_\bb(M;\Lambda)$ is the Reshetikhin-Turaev invariant.

In order to make this construction effective, we need an algorithm for computing the forms $\phi_L$ which are defined by
universal property. Such an algorithm, based on an encoding of certain tangles by means of \emph{Hopf diagrams}, is given in
\cite{BV1}. Thus the invariants $\tau_\bb(M;\alpha)$ can be expressed in terms of certain structural morphisms of the coend
$C$. Section~\ref{Sect-invquant} is devoted to these quantum invariants and their computation.

\medskip
Hence, when $\cc$ is a spherical fusion category, we may compute $\tau_{\zz(\cc)}(M;\Lambda)$ provided we
can describe explicitly the structural morphisms of the coend of $\zz(\cc)$. In other words, we need an algebraic
interpretation of the center construction. If~$\cc$ is braided and has a coend $A$ (which is a Hopf algebra), then the category
$\zz(\cc)$ coincides with the category of (right) $A$-modules in $\cc$. However the difficulty here is that we don't want to
assume $\cc$ is braided. To bypass this difficulty, we use the notion of \emph{Hopf monad} introduced in \cite{BV2}.

Hopf monads generalize Hopf algebras in a non-braided setting. In particular, finite-dimensional Hopf algebras and their
different generalizations (Hopf algebras in braided autonomous categories, quantum bialgebroids, etc...) provide examples of
Hopf monads. If fact, any monoidal adjunction between autonomous categories gives rise to a Hopf monad. It turns out that
much of the theory of finite-dimensional Hopf algebras extends to Hopf monads, see \cite{BV2}. In Section~\ref{Sect-monad}, we recall a few results on Hopf monads.

The whole point of introducing Hopf monads here is that they provide an algebraic interpretation of the center construction
\cite{BV3}. If $\cc$ is a \emph{centralizable} autonomous category, meaning that the coend $Z(X)=\int^{Y \in \cc} \ldual{Y}
\otimes X \otimes Y$ exists for any object~$X$ of $\cc$, then $Z$ is a quasitriangular Hopf monad on $\cc$ and the center~$\zz(\cc)$ coincides, as a braided category, with the category of $Z$-modules in $\cc$. In addition, Drinfeld's double
construction extends naturally to Hopf monads. This theory provides a description of the coend of $\zz(\cc)$. In Section
\ref{Sect-thedouble}, we recall a few facts on the double of Hopf monads.

\medskip
In Section \ref{sect-RT}, we apply the above results to spherical fusion categories. Firstly, we obtain a generalization of
M\"uger's theorem on the modularity of the center of a spherical fusion category $\cc$ to the case where $\dim\cc$ is not
necessarily invertible and~$\kk$ is any commutative ring. Denoting by $\{V_i\}_{i \in I}$ a (finite) representative family of scalar
objects of $\cc$, we get:
\begin{equation*}
Z(X)=\bigoplus_{i \in I} \ldual{V}_i \otimes X \otimes V_i,
\end{equation*}
Moreover $\zz(\cc)$ is centralizable and $\dim\zz(\cc)=(\dim\cc)^2$. The underlying object of the coend of
$\zz(\cc)$ is:
\begin{equation*}
C=\bigoplus_{i,j \in I} \ldual{V}_i \otimes \ldual{V}_j \otimes \lldual{V}_i \otimes V_j\,,
\end{equation*}
and all structural morphisms of $C$ (including its integral $\Lambda\co \un \to C$) can be written down explicitly in $\cc$.
Furthermore, $\zz(\cc)$ is always modular. When $\kk$ is an algebraically closed field and $\dim\cc$ is invertible, then $\zz(\cc)$
is a fusion category and so we recover M\"uger's theorem. However, when $\dim\cc$ is not invertible, $\zz(\cc)$ is a non-semisimple
ribbon category. Nevertheless, in this case, the version $\tau_{\zz(\cc)}(M;\Lambda)$ of the Lyubashenko invariant is still
defined and computable in terms of $\cc$.

\section{Conventions and notations}

\subsection{Autonomous categories}
Monoidal categories are assumed to be strict.

Recall that a \emph{duality} in a monoidal category $(\cc,\otimes,\un)$ is a quadruple $(X,Y,e,d)$, where $X$, $Y$ are
objects of $\cc$, $e\co X \otimes Y \to \un$ (the \emph{evaluation}) and $c \co \un \to Y \otimes X$ (the
\emph{coevaluation}) are morphisms in $\cc$, such that:
\begin{equation*}
(e \otimes \id_X)(\id_X \otimes c)=\id_X \quad \text{and} \quad (\id_Y \otimes e)(c \otimes \id_Y)=\id_Y.
\end{equation*}
Then $(X,e,c)$ is a \emph{left dual} of $Y$, and $(Y,e,c)$ is a \emph{right dual} of $X$.

A \emph{left autonomous} category is a monoidal category for which every object $X$ admits a left dual
$(\ldual{X},\lev_X,\lcoev_X)$. Likewise, a \emph{right autonomous} category is a monoidal category for which every object
$X$ admits a right dual $(\rdual{X},\rev_X,\rcoev_X)$.

An \emph{autonomous category} is a monoidal category which is left and right autonomous. Note that in an autonomous
category, there are canonical isomorphisms:
\begin{align*}
&\ldual{(\rdual{X})}\cong X, && \ldual{(X \otimes Y)} \cong \ldual{Y} \otimes \ldual{X}, && \ldual{\un} \cong \un, \\
&\rdual{(\ldual{X})}\cong X, && \rdual{(X \otimes Y)} \cong \rdual{Y} \otimes \rdual{X}, &&  \rdual{\un}\cong \un.
\end{align*}
Subsequently, in formulae, we will often abstain (by abuse) from writing down these isomorphisms.

\subsection{Sovereign categories}
A \emph{sovereign category} is a left autonomous category endowed with a strong monoidal natural transformation $\phi_X\co X
\to \lldual{X}$. Such a transformation is then an isomorphism. A sovereign category is actually autonomous. Furthermore, in a sovereign category $\cc$, one can define the \emph{left and right traces} of an
endomorphism $f\co X \to X$ as:
\begin{align*}
&\tr_l(f)=\lev_X(\id_{\ldual{X}} \otimes f\phi_X^{-1}) \lcoev_{\ldual{X}} \in \End_\cc(\un), \\
&\tr_r(f)=\rev_X( f\phi_{\rrdual{X}} \otimes \id_{\rdual{X}}) \rcoev_{\rdual{X}} \in \End_\cc(\un),
\end{align*}
and the \emph{left and right dimensions} of an object $X$ as $\dim_l(X)=\tr_l(\id_X) $ and $\dim_r(X)=\tr_r(\id_X)$. We have $\dim_r(X)=\dim_l(\ldual{X})$.

\subsection{Braided categories}
A \emph{braided} category is a monoidal category endowed with a \emph{braiding}, that is, a natural  isomorphism $\tau_{X,Y} \co X \otimes Y \to Y \otimes X$ satisfying: $\tau_{X,Y\otimes Z}=(\id_Y \otimes \tau_{X,Z})(\tau_{X,Y} \otimes \id_Z)$ and $\tau_{X
\otimes Y,Z}=(\tau_{X,Z} \otimes \id_Y)(\id_X \otimes \tau_{Y,Z})$.

\subsection{Ribbon categories}
A \emph{twist} on a braided category $\bb$ is a natural isomorphism $\theta_X \co X \to X$
satisfying: $\theta_{X\otimes Y}=(\theta_X \otimes \theta_Y)\tau_{Y,X}\tau_{X,Y}$. If $\bb$ is braided and autonomous, a
twist $\theta$ on $\bb$ is \emph{self-dual} if $\ldual{(}\theta_X)=\theta_{\ldual{X}}$ (or, equivalently,
$(\theta_X\rdual{)}=\theta_{\rdual{X}}$).

A \emph{ribbon category} is a braided autonomous category endowed with a self-dual twist. A ribbon category is naturally equipped with a sovereign structure such that the left and right traces coincide.

\subsection{Coends}\label{sect-coend}
Let $\cc$, $\dd$ be categories and $F\co \cc^\opp \times \cc \to \dd$ be a functor.

A \emph{dinatural transformation} from the functor $F$ to an object $D$ of $\dd$  is family $d=\{d_X \co F(X,X) \to D\}_{X \in
\Ob(\cc)}$ of morphisms in $\dd$ satisfying the dinaturality condition:
\begin{equation*}
d_X F(f, \id_X)=d_Y F(\id_Y,f)
\end{equation*}
for every morphism $f\co X \to Y$ in $\cc$.

A \emph{coend} of $F$  consists of an object $C$ of $\dd$ and a dinatural
transformation $i$ from $F$ to $C$ which is \emph{universal}, that is, for every dinatural transformation $d$ from $F$ to an object $D$ of $\dd$, there exists a unique morphism $\phi\co C \to D$ such that $d_X=\phi \circ i_X$.

If $F$ admits a coend $(C,i)$, then it is unique (up to unique isomorphism) and one denotes $C=\int^{X\in \cc} F(X,X)$. See \cite{ML1} for details.

\subsection{Coends of autonomous categories}\label{sect-coend-cat}
Let $\cc$ be an autonomous category. If it exists, the coend $C=\int^{X \in \cc} \ldual{X} \otimes X$ of the functor $F\co
\cc^\opp \times \cc \to \cc$ defined by $F(X,Y)=\ldual{X}\otimes Y$ is called the \emph{coend of~$\cc$}. The object $C$ is
then a coalgebra in $\cc$ which coacts universally on the objects of $\cc$ via the the (right) coaction:
\begin{equation*}
\delta_X=(\id_X \otimes i_X)(\lcoev_X \otimes \id_X) \co X \to X \otimes C,
\end{equation*}
where $i_Y \co \ldual{Y} \otimes Y \to C$ is the universal dinatural transformation.

Furthermore, when $\cc$ is braided, $C$ is a Hopf algebra in $\cc$ (see \cite{Maj1,Lyu1}).

\subsection{Dimension of sovereign categories}
Let $\cc$ be a sovereign category which admits a coend. The \emph{left} and \emph{right dimensions of $\cc$} are defined respectively as the left and right dimensions of its coend. These dimensions are actually independent of the choice of sovereign structure on $\cc$. If they coincide (for instance when $\cc$ is a ribbon category or $\cc$ is a fusion category), they are called the \emph{dimension of $\cc$} and denoted $\dim \cc$.

\subsection{Fusion categories}\label{sect-fusion} A \emph{fusion category} over a commutative ring $\kk$ is
a \kt linear autonomous category $\cc$ endowed with a finite family $\{V_i\}_{i \in I}$ of objects of $\cc$ satisfying:
\begin{itemize}
\item $\Hom_\cc(V_i,V_j)=\delta_{i,j} \, \kk$ for all $i,j \in I$;
\item each object of $\cc$ is a finite direct sum of objects of $\{V_i\}_{i \in I}$;
\item $\un$ is isomorphic to some $V_0$ with $0 \in I$.
\end{itemize}

An object $X$ of $\cc$ is \emph{scalar} if  $\End(X)=\kk$. The family $\{V_i\}_{i \in I}$ is a representative family of scalar
objects of $\cc$.
Left and right dualities in $\cc$ preserve scalar objects, and so induce bijections $i \mapsto \ldual{i}$ and $i \mapsto \rdual{i}$ of~$I$ such
that $\ldual{(V_i)} \cong V_{\ldual{i}}$ and $\rdual{(V_i)} \cong V_{\rdual{i}}$. Note that $\ldual{0}=0=\rdual{0}$.

Let $\cc$ be a fusion category. The $\Hom$ spaces in $\cc$ are
free $\kk$-modules of finite type. The \emph{multiplicity} of $i\in I$ in an objet $X$ of $\cc$ is defined as:
\begin{equation*}
N^i_X=\rank_\kk \,\Hom_\cc(V_i,X)=\rank_\kk \,\Hom_\cc(X,V_i).
\end{equation*}
Note there exist morphisms $(p_X^{i,\alpha} \co X \to V_i)_{1 \leq \alpha \leq N^i_X}$ and $(q_X^{i,\alpha}\co V_i \to X)_{1
\leq \alpha \leq N^i_X}$ such that:
\begin{equation*}
\id_X=\!\!\!\!\!\sum_{\substack{i \in I \\ 1 \leq \alpha \leq N^i_X}} \!\!\!q_X^{i,\alpha}p_X^{i,\alpha} \quad \text{and}
\quad p_X^{i,\alpha}q_X^{j,\beta}=\delta_{i,j}\delta_{\alpha,\beta} \, \id_{V_i}.
\end{equation*}

A fusion category $\cc$ admits a coend  $C=\bigoplus_{i \in I} \ldual{V_i} \otimes V_i$ with universal dinatural transformation
given by:
\begin{equation*}
i_X=\sum_{\substack{i\in I \\ 1 \leq \alpha \leq N_X^i}} \ldual{q}_X^{i,\alpha}\otimes p_X^{i,\alpha}.
\end{equation*}
Since $\dim_l(C)=\dim_r(C)$,  the dimension of a sovereign fusion category $\cc$ is:
\begin{equation*}
\dim \cc=\sum_{i \in I} \dim_l(V_i) \dim_r(V_i) \in \kk.
\end{equation*}

In a sovereign fusion category $\cc$, the dimensions $\dim_l(V_i)$ and $\dim_r(V_i)$ of the scalar objects are invertible. However $\dim \cc$ may be not invertible.

A fusion category $\cc$ is \emph{spherical} if it is sovereign and the left and right traces of endomorphisms in~$\cc$
coincide. This last condition is equivalent to the equality of left and right dimensions of the scalar objects $V_i$ for
$i\in I$. In a spherical category, the left (and right) dimension of an object $X$ is denoted $\dim(X)$.

\section{Quantum invariants and Hopf diagrams}\label{Sect-invquant}
In this section, we review a general construction of quantum invariants (of Reshetikhin-Turaev type) and a method for computing them via Hopf diagrams.

\subsection{Constructing quantum invariants}\label{sect-construct}
Let $\bb$ be a ribbon autonomous category ($\bb$ is not necessarily linear). Assume that  $\bb$ admits a coend:
\begin{equation*}
C=\int^{Y \in \bb} \ldual{Y} \otimes Y,
\end{equation*}
with universal coaction $\delta_Y \co Y \to Y \otimes C$ (see Section~\ref{sect-coend-cat}). In particular, using the
general theory of coends, we have the following universal property: for any natural transformation
$\xi=\{\xi_{Y_1,\dots,Y_n} \co Y_1 \otimes \cdots \otimes Y_n \to Y_1 \otimes \cdots \otimes Y_n \otimes M\}_{Y_1,\dots,Y_n
\in \Ob(\bb)}$, where $M$ is an object of $\bb$, there exists a unique morphism $r\co C^{\otimes n} \to M$ such that:
\begin{equation*}
 \psfrag{X}[Bl][Bl]{\scalebox{.8}{$Y_1$}}
 \psfrag{Y}[Bl][Bl]{\scalebox{.8}{$Y_n$}}
 \psfrag{M}[Bl][Bl]{\scalebox{.8}{$M$}}
 \psfrag{H}[cc][cc]{\scalebox{.8}{$\delta_{Y_1}$}}
 \psfrag{K}[cc][cc]{\scalebox{.8}{$\delta_{Y_n}$}}
 \psfrag{C}[cc][cc]{\scalebox{.8}{$C$}}
 \psfrag{B}[cc][cc]{\scalebox{1}{$r$}}
 \psfrag{R}[cc][cc]{\scalebox{1}{$\xi_{Y_1,\dots,Y_n}$}}
 \scaleraisedraw{.5}{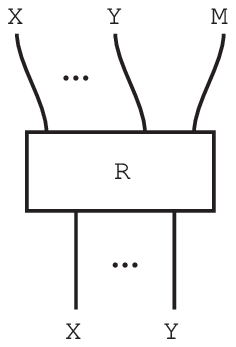} \quad = \quad \scaleraisedraw{.5}{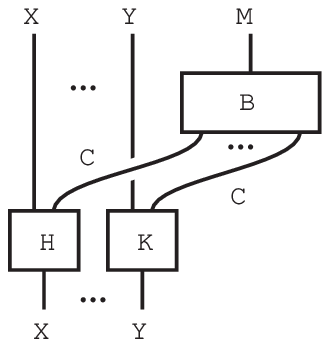}
\end{equation*}
for all objects $Y_1, \dots , Y_n$ of $\bb$.

Now let $T$ be a ribbon $n$-string link with $n$ a non-negative integer. Recall $T$ is a ribbon $(n,n)$-tangle consisting of
$n$ arc components, without any closed component, such that the $k$th arc ($1\leq k \leq n$) joins the $k$th bottom endpoint
to the $k$th top endpoint. We orient $T$ from bottom to top. By virtue of the universality of the category of colored ribbon
tangles, coloring the $n$ components of $T$ with objects $Y_1, \dots, Y_n$ of $\bb$ yields a morphism $T_{Y_1, \cdots, Y_n}
\co Y_1 \otimes \cdots \otimes Y_n \to Y_1 \otimes \cdots \otimes Y_n$, that is,
\begin{equation*}
T_{Y_1, \cdots,Y_n}= \psfrag{X}[cl][cl]{\scalebox{.8}{$Y_1$}} \psfrag{Y}[cl][cl]{\scalebox{.8}{$Y_n$}}
\psfrag{H}[cc][cc]{\scalebox{.8}{$\delta_{Y_1}$}} \psfrag{K}[cc][cc]{\scalebox{.8}{$\delta_{Y_n}$}}
 \psfrag{C}[cc][cc]{\scalebox{.8}{$C$}}  \psfrag{B}[cc][cc]{\scalebox{.8}{$\phi_T$}}
 \scaleraisedraw{.5}{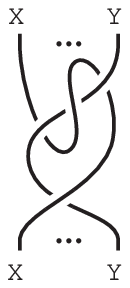}\, .
\end{equation*}
Moreover $T_{Y_1, \cdots, Y_n}$ is natural in each variable $Y_k$ and so, by universality of the coaction of the coend $C$,
there exists a unique morphism:
\begin{equation*}
\phi_T \co C^{\otimes n} \to \un
\end{equation*}
such that:
\begin{equation*}
T_{Y_1, \cdots,Y_n}= \psfrag{X}[cl][cl]{\scalebox{.8}{$Y_1$}} \psfrag{Y}[cl][cl]{\scalebox{.8}{$Y_n$}}
\psfrag{H}[cc][cc]{\scalebox{.8}{$\delta_{Y_1}$}} \psfrag{K}[cc][cc]{\scalebox{.8}{$\delta_{Y_n}$}}
 \psfrag{C}[cc][cc]{\scalebox{.8}{$C$}}  \psfrag{B}[cc][cc]{\scalebox{.8}{$\phi_T$}}
 \;\scaleraisedraw{.5}{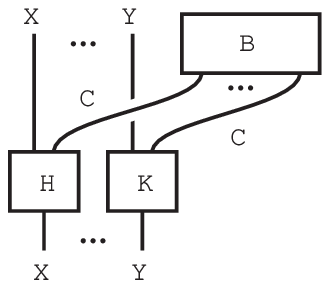} \, .
\end{equation*}

\medskip

Two natural questions arise in this context:
\begin{itemize}
\item How to evaluate the forms $\phi_T$ to get invariants of framed links\footnote{A framed link with $n$ components is
always the closure of some ribbon $n$-string link.} and, further, of 3-manifolds?
\item How to compute the forms $\phi_T$ which are defined by universal property?
\end{itemize}
We address the first question in Section~\ref{sect-kirby} and the second one in Section~\ref{sect-Hdiag}.

\subsection{Kirby elements and quantum invariants}\label{sect-kirby}
As in the previous section, let~$\bb$ be a ribbon autonomous category with a coend $C$. In this setting,
$\Bbbk=\End_\cc(\un)$ is a commutative monoid.

Let $L$ be a framed link in $S^3$ with $n$ components. There always exists a (non-unique) ribbon $n$-string link $T$ such
that $L$ is isotopic to the closure of $T$. For $\alpha \in \Hom_\cc(\un,C)$, set
\begin{equation*}
\tau_\bb(L;\alpha)=\phi_{T} \circ \alpha^{\otimes n} \in \Bbbk,
\end{equation*}
where $\phi_{T}\co C^{\otimes n} \to \un$ is defined as above.

Following~\cite{Vir}, by a \emph{Kirby element} of $\bb$, we mean a morphism $\alpha \in \Hom_\bb(\un,C)$ such that, for any
framed link $L$, $\tau_\bb(L;\alpha)$ is well-defined and invariant under isotopies and 2-handle slides of $L$. A Kirby
element $\alpha$ of $\bb$ is said to be \emph{normalizable} if $\tau_\bb(\bigcirc^{+1};\alpha)$ and
$\tau_\bb(\bigcirc^{-1};\alpha)$ are invertible in $\Bbbk$, where $\bigcirc^{\pm 1}$ denotes the unknot with framing~$\pm
1$.

By universality of the coaction $\delta$ of $C$ on objects of $\cc$, we see that the twist $\theta_Y\co Y \to Y$ of $\bb$
and its inverse lead to morphisms $\theta^\pm_C \co C \to \un$ such that:
\begin{equation*}
\theta_Y^{\pm 1}=(\id_Y \otimes \theta^\pm_C)\delta_Y.
\end{equation*}a
If $\alpha$ is a Kirby element of $\bb$, we have: $\tau_\bb(\bigcirc^{\pm 1};\alpha)=\theta^\pm_C \alpha$, so that $\alpha$
is normalizable if and only if $\theta^\pm_C \alpha$ are invertible in $\Bbbk$.

Recall (see \cite{Li}) that every (closed, connected, oriented) 3-manifold can be obtained from $S^3$ by surgery along a
framed link $L \subset S^3$. For any framed link $L$ in $S^3$, we will denote by $M_L$ the 3-manifold obtained from $S^3$ by
surgery along $L$, by $ n_L$ the number of components of $L$, and by $b_-(L)$ the number of negative eigenvalues of the
linking matrix of $L$.

An immediate consequence of the Kirby theorem \cite{Ki} is that if $\alpha$ is a normalizable Kirby element of $\bb$, then:
\begin{equation*}
\tau_\bb(M_L;\alpha)=(\theta^+_C \alpha)^{b_-(L)-n_L}\,  (\theta^-_C \alpha)^{-b_-(L)} \; \tau_\bb(L;\alpha)
\end{equation*}
is an invariant of 3-manifolds. Furthermore these invariants are multiplicative under the connected sum of 3\ti manifolds:
$\tau_\bb(M\# M';\alpha)=\tau_\bb(M;\alpha)\,\tau_\bb(M';\alpha)$.

Note that if $\alpha$ is a normalizable Kirby element and $k$ is an automorphism of $\un$,
then $k\alpha$ is also a normalizable Kirby element. The normalization of the invariant $\tau_\bb(M;\alpha)$ has been chosen
so that $\tau_\bb(M;k\alpha)=\tau_\bb(M;\alpha)$.\\

The question is now: how to determine the (normalizable) Kirby element of $\bb$? A partial answer was given in \cite{Vir}.
Denoting by $m_C$, $\Delta_C$, and $S_C$ respectively the product, coproduct, and antipode of the Hopf algebra $C$, we have:
\begin{thm}[{\cite[Theorem~2.5]{Vir}}]\label{thmkirbel}
Any morphism $\alpha\co  \un \to C$ in $\bb$ such that:
\begin{equation*}
S_C \alpha=\alpha \quad \text{and} \quad (m_C \otimes \id_C)(\id_C \otimes \Delta_C)(\alpha \otimes \alpha)=\alpha \otimes
\alpha
\end{equation*}
is a Kirby element of $\bb$.
\end{thm}
For instance, the unit $u_C$ of $C$ is a normalizable Kirby element (its associated invariant is the trivial one).

A more interesting example of a a Kirby element is an $S$-invariant integral $\Lambda$ of $C$, that is, a morphism $\Lambda
\co \un \to C$ such that $S_C(\Lambda)=\Lambda$ and $m_C(\Lambda \otimes \id_C)=\Lambda\, \varepsilon_C=m_C(\id_C \otimes
\Lambda)$, where $\varepsilon_C$ is the counit of $C$. For the existence of such integrals,  we refer to \cite{BKLT}. If
$\Lambda$ is normalizable, then the associated invariant is the Lyubashenko's one \cite{Lyu2}, up to a different
normalization.

Note that other Kirby elements exist in general (see \cite{Vir}).

\begin{rem}\label{rem-invRT}
Assume $\bb$ is a modular category in the sense of \cite{Tur2}, that is, a ribbon fusion category with invertible
$S$-matrix. Let $\{V_i\}_{i\in I}$ be a representative family of simple objects of $\bb$. Then $\bb$ admits a coend
$C=\bigoplus_{i \in I} \ldual{V}_i \otimes V_i$. Let $\phi_X\co X \to \lldual{X}$ be the sovereign structure of $\bb$ and
set:
\begin{equation*}
\Lambda=\sum_{i \in I} \dim(V_i) \, (\id_{\ldual{V}_i} \otimes \phi_i^{-1}) \lcoev_{V_i} \co \un \to C,
\end{equation*}
Then $\Lambda$ is a $S_C$-invariant integral of $C$. Furthermore it is normalizable and its associated invariant is the
Reshetikhin-Turaev one~\cite{Tur2}, up to a different normalization. More precisely, assuming $\dim \bb=\sum_{i \in I}
\dim(V_i)^2$ has a square root $D$ in $\Bbbk$ (which is then invertible in this context), setting $\Delta_-=\theta^-_C
\Lambda$, and denoting by $b_1(M)$ the first Betti number of $M$, we have:
\begin{equation*}
\mathrm{RT}_\bb(M) = D^{-1} \Bigl(\frac{D}{\Delta_-}\Bigr)^{b_1(M)} \, \tau_\bb(M;\Lambda).
\end{equation*}
We will see in Section~\ref{sect-RT} that, unlike $\mathrm{RT}_\bb(M)$, $\tau_\bb(M;\Lambda)$ may be still defined for
ribbon categories $\bb$ with $\dim \bb=0$.
\end{rem}

\subsection{Hopf diagrams}\label{sect-Hdiag}
For a precise treatment of the theory of Hopf diagrams, we refer to \cite{BV1}. Note that Habiro, shortly after us, had
similar results in \cite{Hab1}.

\begin{figure}[t]
   \begin{center}
       $\Delta=\,$\scaleraisedraw{.403}{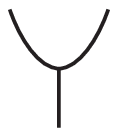}\,,  \qquad $\varepsilon=\;$\scaleraisedraw{.403}{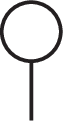}\,, \qquad $\omega_+=\;$\scaleraisedraw{.403}{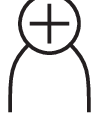}\,, \qquad $\omega_-=\;$\scaleraisedraw{.403}{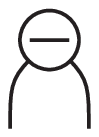}\,, \\[.9em]
       $\theta_+=\;$\scaleraisedraw{.403}{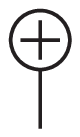}\,,
         \qquad $\theta_-=\;$\scaleraisedraw{.403}{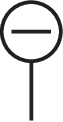}\,,  \qquad $S=\;$\scaleraisedraw{.403}{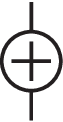}\,,  \qquad $S^{-1}=\;$\scaleraisedraw{.403}{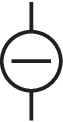}\,, \\[.9em]
          $\tau=\,$\scaleraisedraw{.403}{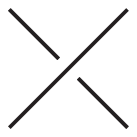}\,,  \qquad \quad $\tau^{-1}=\,$\scaleraisedraw{.403}{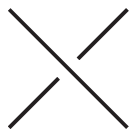}\,.
   \end{center}
     \caption{Generators of Hopf diagrams}
     \label{gensHdiag}
\end{figure}

\begin{figure}[t]
   \begin{center}
       \subfigure[A Hopf diagram with 1 input]{\phantom{XXXXXX}\scaledraw{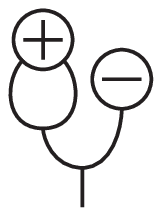}\phantom{XXXXXX}} \quad\;
       \subfigure[A Hopf diagram with 2 inputs]{\phantom{XXXXX}\scaledraw{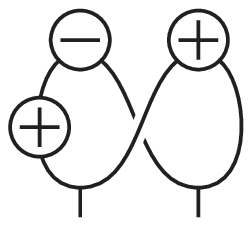}\phantom{XXXXX}}
   \end{center}
   \caption{Examples of Hopf diagrams}
     \label{exsHdiag}
\end{figure}

\begin{figure}[t]
   \begin{center}
       \scaleraisedraw{.403}{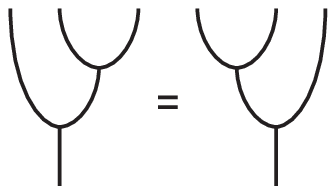}\quad \text{denoted} \quad \scaleraisedraw{.403}{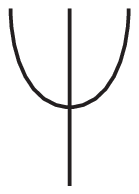}\, , \qquad
       \scaleraisedraw{.403}{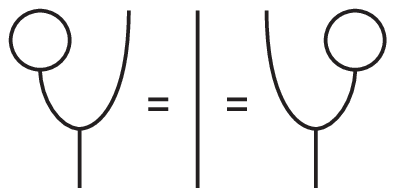}, \\[1.1em]
       \scaleraisedraw{.403}{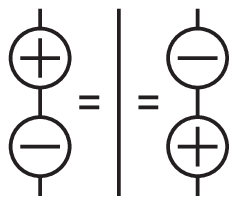}\; , \qquad
       \scaleraisedraw{.403}{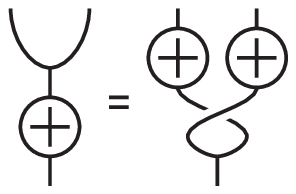} , \qquad
       \scaleraisedraw{.403}{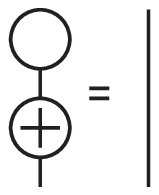}\; , \qquad
       \scaleraisedraw{.403}{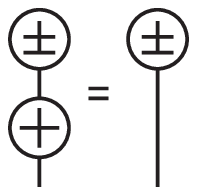} , \\[1.1em]
       \scaleraisedraw{.403}{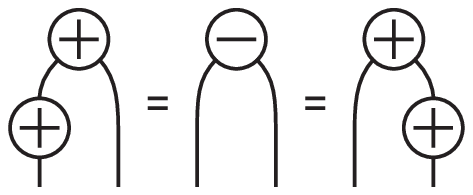}\; , \qquad
       \scaleraisedraw{.403}{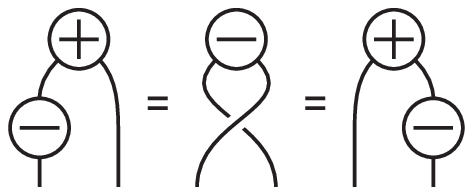}\; , \\[1.1em]
       \psfrag{,}[Bl][Bl]{,}
       \psfrag{.}[Bl][Bl]{.}
       \scaledraw{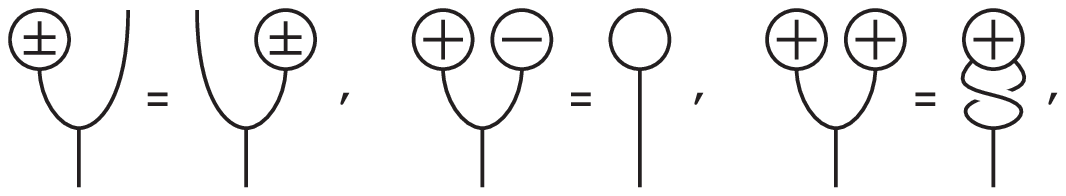}\\[1em]
       \scaledraw{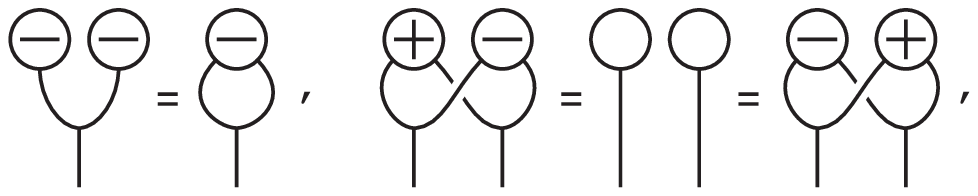} \\[1em]
       \scaledraw{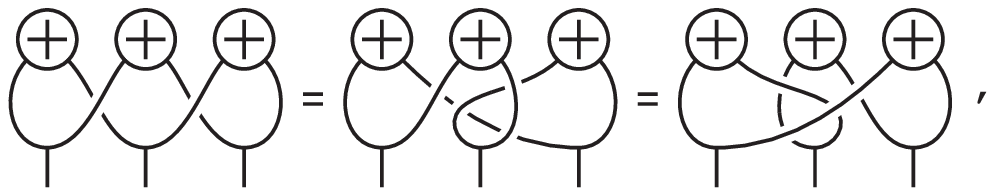} \\[1em]
       \scaledraw{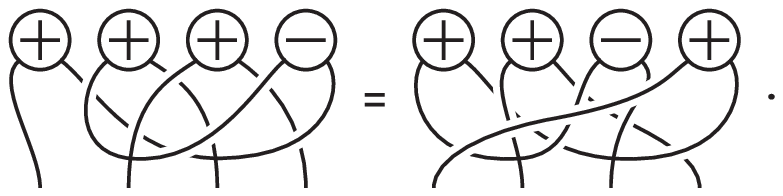}
   \end{center}
     \caption{Relations on Hopf diagrams}
     \label{relsHdiag}
\end{figure}

Briefly speaking, a \emph{Hopf diagram} is a planar diagram, with inputs but no output, obtained by stacking the generators
of Figure~\ref{gensHdiag} (diagrams are read from bottom to top). Examples of Hopf diagrams with~1 and 2 inputs are depicted
in Figure~\ref{exsHdiag}. Hopf diagrams are submitted to the relations of Figure~\ref{relsHdiag} (plus relations expressing
that $\tau$ is an invertible QYBE solution which is natural with  respect to the other generators). In particular, the
relations of Figure~\ref{relsHdiag} say that $\Delta$ behaves as a coproduct with counit $\varepsilon$, $S$ behaves as an
antipode, $\omega_\pm$ behaves as a Hopf pairing, and $\theta_\pm$ behaves as a twist form. The last two relations of
Figure~\ref{relsHdiag} are nothing but the Markov relations for pure braids.

\begin{figure}[h]
   \begin{center}
       \psfrag{F}[cc][cc]{$D\star D'$}
       \psfrag{D}[cc][cc]{$D$}
       \psfrag{E}[cc][cc]{$D'$}
       \psfrag{o}[cc][cc]{$\circ$}
       \scaledraw{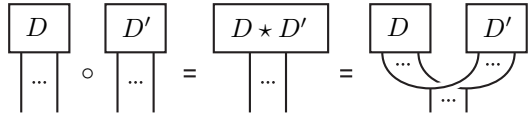}
   \end{center}
     \caption{Composition of Hopf diagrams}
     \label{convoprod}
\end{figure}

Hopf diagrams with the same number of inputs can be composed using the \emph{convolution product} $\star$ defined in
Figure~\ref{convoprod}. This leads to the category $\diag$ of Hopf diagrams. Objects of $\diag$ are the non-negative
integers. For two non-negative integers $m$ and $n$, the set $\Hom_{\diag}(m,n)$ of morphisms from $m$ to $n$ is the empty
set if $m \neq n$ and is the set of Hopf diagrams with $m$ inputs (up to their relations) if $m=n$. The composition is the
convolution product and the identity of $n$ is the Hopf diagram obtained by juxtaposing $n$ copies of $\varepsilon$.

The category $\diag$ is a monoidal category: $m \otimes n=m+n$ on objects and the monoidal product $D \otimes D'$ of two
Hopf diagrams $D$ and $D'$ is the Hopf diagram obtained by juxtaposing $D$ on the left of $D'$.

Let us denote by $\rstl$ the category of ribbon string links. The objects of $\rstl$ are the non-negative integers. For two
non-negative integers $m$ and $n$, the set of morphisms from $m$ to $n$ is
\begin{equation*}
\Hom_{\rstl}(m,n)=\begin{cases} \emptyset & \text{if $m\neq n$,}\\ \rstl_n & \text{if $m=n$,} \end{cases}
\end{equation*}
where $\rstl_n$ denotes the set of (isotopy classes) of ribbon $n$-string links. The composition $T' \circ T$ of two ribbon
$n$-string links is given by stacking $T'$ on the top of $T$ (i.e., with ascending convention). Identities are the trivial
string links. Note that the category $\rstl$ is a monoidal category: $m \otimes n=m+n$ on objects and the monoidal product
$T \otimes T'$ of two ribbon string links $T$ and $T'$ is the ribbon string link obtained by juxtaposing $T$ on the left of
$T'$.

\begin{figure}[t]
   \begin{center}
      \psfrag{,}[Bl][Bl]{,}
      \psfrag{.}[Bl][Bl]{.}
      \scaledraw{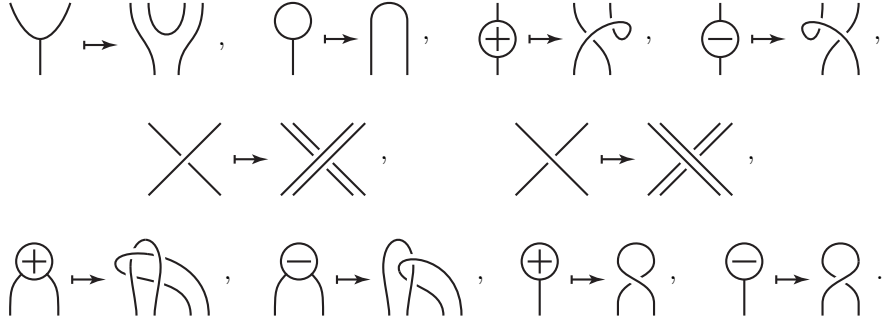}
   \end{center}
     \caption{Rules for transforming Hopf diagrams to tangles}
     \label{Hopf2Hand}
\end{figure}

\begin{figure}[t]
   \begin{center}
$D=$ \scaleraisedraw{.403}{2imp} \hspace*{.1pt} $\rightsquigarrow$  $h_D=$   \scaleraisedraw{.403}{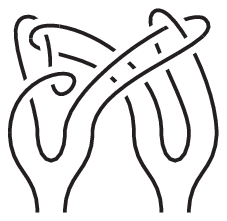} \hspace*{.1pt}
$\rightsquigarrow$ $\Phi(D)=$\, \scaleraisedraw{.403}{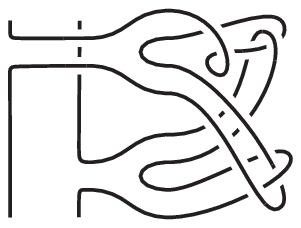} $\sim$ \,\scaleraisedraw{.403}{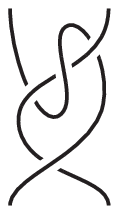}
   \end{center}
     \caption{From Hopf diagrams to ribbon string links}
     \label{exHopf2Hand}
\end{figure}

Hopf diagrams give a `Hopf algebraic' description of ribbon string links. Indeed, any Hopf diagram $D$ with $n$ inputs gives
rise to a ribbon $n$-string link $\Phi(D)$ in the following way: using the rules of Figure~\ref{Hopf2Hand}, we obtain a
ribbon $n$-handle\footnote{Ribbon handles are called bottom tangles in \cite{Hab1}.} $h_D$, that is, a ribbon
$(2n,0)$-tangle consisting of $n$ arc components, without any closed component, such that the $k$\trait th arc joins the
$(2k-1)$\trait to the $2k$\trait th bottom endpoints. Then, by rotating $h_D$, we get a ribbon $n$-string link $\Phi(D)$:
\begin{equation*}
\text{$D$ Hopf diagram} \;\; \rightsquigarrow  \; \;   \psfrag{a}[cc][cc]{$h_D$} \scaleraisedraw{.403}{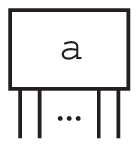} \;\;
\rightsquigarrow \;\; \Phi(D)=\, \scaleraisedraw{.403}{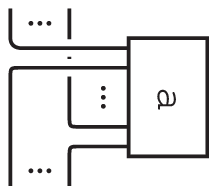} \; .
\end{equation*}
An example of this procedure is depicted in Figure~\ref{exHopf2Hand}.

This leads to a functor $\Phi\co \diag \to \rstl$ defined on objects by $n \mapsto \Phi(n)=n$ and on morphisms by $D \mapsto
\Phi(D)$.
\begin{thm}[{\cite[Theorem 4.5]{BV1}}]
$\Phi\co \diag \to \rstl$ is a well-defined monoidal functor and there exists (constructive proof) a monoidal functor
$\Psi\co \rstl \to \diag$ which satisfies $\Phi\circ\Psi=1_{\rstl}$.
\end{thm}
Note that by `constructing proof' we mean there is an explicit algorithm that associates to a ribbon string $T$ a Hopf
diagram $\Psi(T)$ such that $\Phi\bigl(\Psi(T) \bigr)=T$ (see~\cite{BV1}). The key point is that such a functor $\Psi$
exists thanks to the relations we put on Hopf diagrams.\\

Let now $\bb$ be a ribbon autonomous category which admits a coend $C$. Let us answer to the second question of
Section~\ref{sect-construct}: given a ribbon $n$-string link $T$, how to compute the morphism $\phi_T\co C^{\otimes n}\to
\un$ which is defined by universal property? Recall $C$ is a Hopf algebra in $\bb$ and denote its coproduct, counit, and
antipode by $\Delta_C$, $\varepsilon_C$, and $S_C$ respectively. The twist (and its inverse) of $\bb$ is encoded by
morphisms $\theta^\pm_C \to C \to \un$ (see Section~\ref{sect-kirby}). Furthermore, we can defined a Hopf pairing $\omega_C
\co C \otimes C \to \un$ via:
\begin{align*}
&\omega_C (i_X \otimes i_Y)=(\lev_X \otimes \lev_Y)(\id_{\ldual{X}} \otimes \tau_{\ldual{Y},X}\tau_{X,\ldual{Y}} \otimes
\id_{\ldual{Y}}),
\end{align*}
where $\tau$ is the braiding of $\bb$ and $i_Y \co \ldual{Y} \otimes Y \to C$ is the universal dinatural transformation of
the coend $C$. Finally, we set $\omega_C^+=\omega_C(S_C^{-1} \otimes \id_C)$ and $\omega_C^-=\omega_C$.

\begin{thm}[{\cite[Theorem 5.1]{BV1}}]
Let $T$ be ribbon $n$-string link. Let $D$ be any Hopf diagram (with $n$ entries) which encodes $T$, that is, such that
$\Phi(D)=T$ (recall there is an algorithm producing such a Hopf diagram). Then the morphism $\phi_T\co C^{\otimes n} \to
\un$ defined by $T$ is given by replacing in $D$ the generators $\Delta$, $\varepsilon$, $\omega_\pm$, $\theta_\pm$, $S^{\pm
1}$, and $\tau^{\pm 1}$ (see Figure~\ref{gensHdiag}) by the morphisms $\Delta_C$, $\varepsilon_C$, $\omega^\pm_C$,
$\theta^\pm_C$, $S^{\pm 1}_C$, and $\tau^{\pm 1}_{C,C}$ respectively.
\end{thm}

Remark that the product and unit of the coend $C$ are not needed to represent Hopf diagrams.\\

Let us summarize the above universal construction of  quantum invariants, starting from a ribbon category $\bb$ which admits
a coend $C$. Pick a normalizable Kirby element $\alpha$ of $\bb$ (for example as in Theorem~\ref{thmkirbel}). Recall it
gives rise to the invariant $\tau_\bb(M,\alpha)$ of 3\ti manifolds. Let $M$ be a 3\ti manifold. Present $M$ by surgery along
a framed link $L$, which can be viewed as the closure of a ribbon $n$-string link $T$ where $n$ is the number of components
of $L$. Encode the string link $T$ by a Hopf diagram~$D$:
\begin{equation*}
M\simeq S^3_L, \quad L\, \sim \psfrag{T}[c][c]{\scalebox{.72}{$T$}}\sdraw{.38}{1}{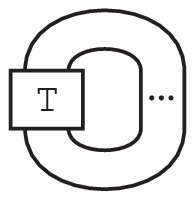} \quad \text{with} \quad
T=\,\sdraw{.38}{1}{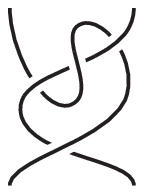} \; \scalebox{-1.1}[1.1]{$\mapsto$ \,}\quad D=\sdraw{.38}{1}{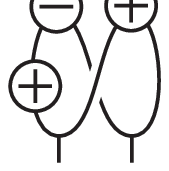}
\end{equation*}
The morphism $\phi_T\co C^{\otimes n} \to \un$ associated to $T$ can be computed by replacing the generators of $D$ by the
corresponding structural morphisms of the coend $C$. Then evaluate $\phi_T$ with the Kirby element $\alpha$ and normalize to
get the invariant:
\begin{equation*}
\tau_\bb(M;\alpha)=\;\psfrag{w}[c][c]{\scalebox{.8}{$\omega^-_C$}} \psfrag{m}[c][c]{\scalebox{.8}{$\omega^+_C$}}
\psfrag{o}[c][c]{\scalebox{.8}{$S_C$}} \psfrag{c}[c][c]{\scalebox{.8}{$\tau_{C,C}$}}
\psfrag{D}[c][c]{\scalebox{.8}{$\Delta_C$}} \psfrag{a}[c][c]{\scalebox{.9}{$\alpha$}}
\psfrag{t}[c][c]{\scalebox{.8}{$\theta^+_C$}} \psfrag{r}[c][c]{\scalebox{.8}{$\theta^-_C$}}
\psfrag{x}[Bl][Bl]{\scalebox{.8}{$b_-(L)-n$}} \psfrag{u}[Bl][Bl]{\scalebox{.8}{$-b_-(L)$}}\scaleraisedraw{.403}{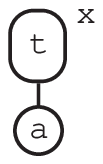}
\qquad \quad\;\;\scaleraisedraw{.403}{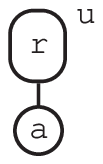} \quad \qquad\scaleraisedraw{.403}{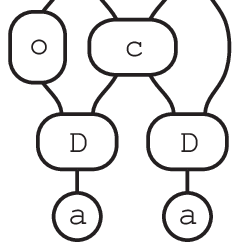}
\end{equation*}

In particular, to compute such quantum invariants defined from the center $\zz(\cc)$ of a autonomous category $\cc$, one
needs to give an explicit description of the structural morphism of the coend of $\zz(\cc)$ in terms of the category $\cc$.
In the next section, we give such a description by using Hopf monads (this was our original motivation for introducing Hopf
monad).

\section{Hopf monads}\label{Sect-monad}

In this section, we review some facts on Hopf monads \cite{BV2}.

\subsection{Monads}
Let $\cc$ be a category. Recall that the category $\End(\cc)$ of endofunctors of $\cc$ is strict monoidal with composition
for monoidal product and identity functor $1_\cc$ for unit object.

A \emph{monad} on $\cc$  (also called a \emph{triple}) is an algebra in $\End(\cc)$, that is, a triple $(T,\mu,\eta)$, where
$T\co \cc \to \cc$ is a functor, $\mu\co T^2 \to T$ and $\eta\co 1_\cc \to T$ are natural transformations, such that:
\begin{equation*}
\mu_X T(\mu_X)=\mu_X\mu_{T(X)} \quad \text{and} \quad \mu_X\eta_{T(X)}=\id_{T(X)}=\mu_X T(\eta_X)
\end{equation*}
for any object $X$ of $\cc$.

\subsection{Modules over a monad} Let $T$ be a monad on a category $\cc$. An \emph{action} of~$T$ on an object $M$
of $\cc$ is a morphism $r\co T(M) \to M$ in $\cc$ such that:
\begin{equation*}
r T(r)= r \mu_M \quad \text{and} \quad r \eta_M= \id_M.
\end{equation*}
The pair $(M,r)$ is then called a \emph{$T$-module in $\cc$}, or just a \emph{$T$-module}\footnote{This is not standard
terminology: pairs $(M,r)$ are usually called $T$-algebras in the literature. However pairs $(M,r)$ are considered here as
the analogues of modules over an algebra, and so the term `algebra' would be awkward in this context.}.

Given two $T$-modules $(M,r)$ and $(N,s)$ in $\cc$, a morphism $f\in \Hom_\cc(M,N)$ is said to be \emph{$T$-linear} if $f
r=s T(f)$. This gives rise to the \emph{category $T\ti\cc$ of $T$-modules}, with composition inherited from $\cc$.

We will denote $U_T\co T\ti\cc \to \cc$ the \emph{forgetful functor of $T$} defined by $U_T(M,r)=M$  for any $T$-module
$(M,r)$ and $U_T(f)=f$ for any $T$-linear morphism~$f$.

\subsection{The philosophy}
Roughly speaking, a monad $T$ on a monoidal category $\cc$ is a bimonad, a Hopf monad, a quasitriangular Hopf monad, or a
ribbon Hopf monad if the category $T\ti\cc$ of $T$-modules is respectively monoidal, autonomous, braided, or ribbon, in such
a way the forgetful functor $U_T \co T\ti\cc \to \cc$ is strict monoidal.

The key point is that these categorical properties of $T\ti\cc$ can be encoded by structural morphisms of $T$. In the next
sections, we briefly give the definitions of these structural morphisms. Their relations with the category $T\ti\cc$ is
summarized in Theorem~\ref{thm-biHopfmon}. For a complete treatment, we refer to \cite{BV2}.

\subsection{Bimonads}
A \emph{bimonad}\footnote{This notion of bimonad coincides exactly with the notion of `Hopf monad' introduced in
\cite{Moer}. However, by analogy with the notions  of bialgebra and Hopf algebra, we prefer to reserve the term `Hopf monad'
for bimonads with antipodes (see Section~\ref{sect-hopfmon})} on a monoidal category $\cc$ is a monad $(T,\mu,\eta)$
on~$\cc$ endowed with a natural transformation $ T_2(X,Y)\co T(X \otimes Y) \to T(X) \otimes T(Y)$ and a morphism $T_0 \co
T(\un) \to \un$ satisfying:
\begin{align*}
& (\id_{T(X)} \otimes T_2(Y,Z))  T_2(X,Y \otimes Z)= (T_2(X,Y) \otimes \id_{T(Z)}) T_2(X \otimes Y, Z);\\
& (\id_{T(X)} \otimes T_0) T_2(X,\un)=\id_{T(X)}=(T_0 \otimes \id_{T(X)}) T_2(\un,X); \\
& T_2(X,Y) \mu_{X \otimes Y}= (\mu_X \otimes \mu_Y) T_2(T(X),T(Y)) T(T_2(X,Y)) ;\\
& T_0 \mu_\un= T_0 T(T_0); \qquad T_2(X,Y) \eta_{X \otimes Y}= (\eta_X \otimes \eta_Y); \qquad  T_0 \eta_\un= \id_\un;
\end{align*}
for all objects $X,Y,Z$ of $\cc$.

\subsection{Antipodes}\label{sect-antipod}
Let $(T,\mu,\eta)$ be a bimonad on a monoidal category $\cc$.

If $\cc$ is left autonomous, then a \emph{left antipode for $T$} is a natural transformation $s^l=\{s^l_X\co T(\ldual{T(X)})
\to \ldual{X}\}_{X \in \Ob(\cc)}$ satisfying:
\begin{align*}
& T_0 T(\lev_X)T(\ldual{\eta_X} \otimes \id_X)=\lev_{T(X)}(s^l_{T(X)}T(\ldual{\mu}_X) \otimes
\id_{T(X)})T_2(\ldual{T(X)},X);  \\
& (\eta_X \otimes \id_{\ldual{X}})\lcoev_X T_0=(\mu_X \otimes s^l_X) T_2(T(X),\ldual{T(X)})T(\lcoev_{T(X)}).
\end{align*}

Likewise, if $\cc$ is right autonomous, then a \emph{right antipode for $T$} is a natural transformation $s^r=\{s^r_X\co
T(\rdual{T(X)}) \to \rdual X\}_{X \in \Ob(\cc)}$ satisfying:
\begin{align*}
& T_0 T(\rev_X)T(\id_X \otimes \eta_X^\vee)=\rev_{T(X)}(\id_{T(X)}
\otimes s^r_{T(X)}T(\mu_X^\vee))T_2(X,\rdual{T(X)}); \\
& (\id_{X^\vee}\otimes \eta_X )\rcoev_X T_0=(s^r_X \otimes\mu_X) T_2(\rdual{T(X)},T(X))T(\rcoev_{T(X)}).
\end{align*}

As in the classical case, left and right antipodes are `anti-(co)multiplicative', see~\cite[Theorem~3.7]{BV2}.

Note that if a left (resp.\@ right) antipode exists, then it is unique. Furthermore, when they exist, the left antipode
$s^l$ and the right antipode $s^r$ are `inverse' to each other in the sense that $ \id_{T(X)}=s^r_{\ldual{T(X)}}
T(\rdual{(s_X^l)})=s^l_{\rdual{T(X)}} T(\ldual{(s_X^r)}) $ for any object~$X$ of $\cc$.

\subsection{Hopf monads}\label{sect-hopfmon}
A \emph{Hopf monad} is a bimonad on an autonomous category which has a left antipode and a right antipode.

Hopf monads generalize Hopf algebras to a non-braided (and non-linear) setting. Furthermore they are much more general: for
example, if $\cc,\dd$ are two autonomous categories and $U\co\dd \to \cc$ is a strong monoidal functor which admits a left
adjoint $F\co \cc \to \dd$, then $T=UF$ is a Hopf monad on $\cc$ (see~\cite[corollary~3.15]{BV2}).

Note that many fundamental results of the theory of Hopf algebras (such as the decomposition of Hopf modules, the existence
of integrals, Maschke's criterium of semisimplicity, etc...)\@ can be generalized to Hopf monads (see \cite{BV2}).

\subsection{Quasitriangular Hopf monads}
Let $T$ be a Hopf monad on an autonomous category~$\cc$. An \emph{\Rt matrix} for $T$ is a natural transformation $R_X,Y\co
X \otimes Y \to T(Y) \otimes T(X)$ satisfying:
\begin{align*}
& (\mu_Y \otimes \mu_X)R_{T(X),T(Y)}T_2(X,Y)=(\mu_Y \otimes \mu_X)T_2(T(Y),T(X))T(R_{X,Y}); \\
\begin{split}
& (\id_{T(Z)}  \otimes T_2(X,Y))R_{X \otimes Y,Z} \\
& \phantom{XXXXXX}=(\mu_Z \otimes \id_{T(X) \otimes T(Y)}) (R_{X,T(Z)} \otimes \id_{T(Y)} )(\id_X \otimes R_{Y,Z}) ;
\end{split}\\
\begin{split}
& (T_2(Y,Z) \otimes \id_{T(X)})R_{X,Y \otimes Z} \\
& \phantom{XXXXXX}=( \id_{T(Y) \otimes T(Z)} \otimes \mu_X) ( \id_{T(Y)} \otimes R_{T(X),Z} ) (R_{X,Y} \otimes \id_Z).
\end{split}
\end{align*}

Note that an \Rt matrix satisfies some QYB equation and is $*$-invertible (where $*$ is some convolution product),
see~\cite[corollary~8.7]{BV2}.

A \emph{quasitriangular Hopf monad} is a Hopf monad equipped with an \Rt matrix.

\subsection{Ribbon Hopf monads} Let $T$ be a quasitriangular Hopf monad $T$ on an autonomous category $\cc$.
A \emph{twist} for $T$ is a central and $*$-invertible natural transformation $\theta_X \co X \to T(X)$ satisfying:
\begin{equation*}
T_2(X,Y) \theta_{X \otimes Y}=(\mu_X \theta_{T(X)} \mu_X \otimes \mu_Y \theta_{T(Y)} \mu_Y)
  R_{T(Y),T(X)} R_{X,Y}.
\end{equation*}
Here central and $*$-invertible means central and invertible in the monoid $\Nat(1_\cc,T)$ of natural transformations from
$1_\cc$ to $T$. This monoid is endowed with the convolution product, defined by:~$(\phi
* \psi)_X=\mu_X \phi_{T(X)} \psi_X =\mu_X T(\psi_X) \phi_X \co X \to T(X)$, and with the unit~$\eta$.

A twist of a quasitriangular Hopf monad on an autonomous category is said to be \emph{self-dual} if it satisfies:
\begin{equation*}
\ldual{\theta_X}=s^l_X\theta_{\ldual{T(X)}} \quad \text{(or, equivalently, $\rdual{\theta}_X=s^r_X\theta_{\rdual{T(X)}}$).}
\end{equation*}

A \emph{ribbon Hopf monad} is a quasitriangular Hopf monad on an autonomous category endowed with a self-dual twist.

\subsection{Relations with modules}
The notions of bimonads, Hopf monads, quasitriangular Hopf monads, or ribbon Hopf monads have a natural interpretation in
terms of the category of modules over the underlying monad. We summarize these properties in the following theorem:
\begin{thm}[\cite{BV2}]\label{thm-biHopfmon}
 \begin{enumerate}
  \renewcommand{\labelenumi}{{\rm (\alph{enumi})}}
\item Let $T$ be a monad on a monoidal category~$\cc$. If~$T$ is a bimonad, then the category $T\ti\cc$ of $T$-modules is
monoidal by setting:
\begin{equation*}
(M,r) \otimes_{T\ti\cc} (N,s)=(M \otimes N, (r \otimes s) T_2(M,N)) \quad \text{and} \quad \un_{T\ti\cc}=(\un,T_0).
\end{equation*}
Moreover this gives a bijective correspondence between bimonad structures for the monad $T$ and  monoidal structures of
$T\ti\cc$ such that the forgetful functor $U_T \co T\ti\cc \to \cc$ is strict monoidal.
  \item Let $T$ be a bimonad on a left autonomous $\cc$. Then $T$ has a left antipode $s^l$ if and only if
     the category $T\ti\cc$ of $T$-modules is
    left autonomous. In terms of a left antipode $s^l$, left duals in
$T\ti\cc$ are given by:
\begin{equation*}
\ldual{(M,r)}=(\ldual{M}, s^l_M T(\ldual{r})),\quad \lev_{(M,r)}=\lev_M, \quad\lcoev_{(M,r)}=\lcoev_M.
\end{equation*}
  \item Let $T$ be a bimonad on a right autonomous $\cc$. Then $T$ has a right antipode $s^l$ if and only if
   the category $T\ti\cc$ of $T$-modules is right autonomous. In terms of a right antipode $s^r$, right duals
in $T\ti\cc$ are given by:
\begin{equation*}
\rdual{(M,r)}=(\rdual{M}, s^r_M T(\rdual{r})),\quad \rev_{(M,r)}=\rev_M, \quad\rcoev_{(M,r)}=\rcoev_M.
\end{equation*}
\item Let $T$ be a bimonad on an autonomous $\cc$. Then $T$ is a Hopf monad if and only if
   the category $T\ti\cc$ of $T$-modules is autonomous.
\item Let $T$ be a bimonad on a monoidal category~$\cc$. Any \Rt matrix $R$ for $T$ yields a braiding $\tau$ on $T\ti\cc$ as
follows:
\begin{equation*}
\tau_{(M,r),(N,s)}=(s \otimes t)R_{M,N}\co (M,r) \otimes (N,s) \to (N,s) \otimes (M,r).
\end{equation*}
This assignment gives a bijection between \Rt matrices for~$T$ and braidings on $T\ti\cc$.
\item Let $T$ be a quasitriangular Hopf monad on an autonomous category~$\cc$. Any twist $\theta$ for $T$ yields a twist $\Theta$
on $T\ti\cc$ as follows:
\begin{equation*}
\Theta_{(M,r)}=r \theta_M\co (M,r) \to (M,r).
\end{equation*}
This assignment gives a bijection between twists for $T$ and twists on $T\ti\cc$. Moreover, in this correspondence, $\theta$
is self-dual (and so $T$ is ribbon) if and only if $\Theta$ is self-dual (and so $T\ti\cc$ is ribbon).
\end{enumerate}
\end{thm}

\section{Quantum double of Hopf monads}\label{Sect-thedouble}
In this section, we review the construction of the double of a Hopf monad and its relations with the center construction
(see \cite{BV3} for details).

\subsection{The center of an monoidal category category}\label{sect-centerDr} Let $\cc$ be a braided category.
Recall that the \emph{center of $\cc$} is the category $\zz(\cc)$ defined as follows: the objects are pairs $(M,\sigma)$,
where $M$ is an object of $\cc$ and $\sigma_Y\co M \otimes Y \to Y \otimes M$ is a natural isomorphism verifying $\sigma_{Y
\otimes Z} =(\id_{Y} \otimes \sigma_Z)(\sigma_Y \otimes \id_{Z})$. A morphism $f \co (M,\sigma) \to (M',\sigma')$ in
$\zz(\cc)$ is a morphism $f\co M \to M'$ in $\cc$ which satisfies $(\id_{Y} \otimes f) \sigma_Y=\sigma'_Y(f \otimes \id_Y)$.
The composition and identities are inherited from that of $\cc$.

The center  $\zz(\cc)$ of $\cc$ is monoidal  with unit object $(\un,\id_M )$ and monoidal product defined by $ (M,\sigma)
\otimes (N,\gamma)=\bigl (M \otimes N,(\sigma \otimes \id_N)(\id_M \otimes \gamma)\bigr)$. Furthermore, if $\cc$ is
autonomous, then so is $\zz(\cc)$.

We define the forgetful functor $\uu\co \zz(\cc) \to \cc$ by $\uu(M,\sigma)=M$ and $\uu(f)=f$. This is a strict monoidal
functor.

\subsection{The double of a Hopf monad}
Let $T$ be a Hopf monad on an autonomous category $\cc$. Assume $T$ is \emph{centralizable}, that is, such that the coend:
\begin{equation*}
Z_T(X)=\int^{Y \in \cc} \ldual{T(Y)} \otimes X \otimes Y
\end{equation*}
exists for every object $X$ of $\cc$. Denote $i_{X,Y} \co \ldual{T(Y)} \otimes X \otimes Y \to Z_T(X)$ the associated
universal dinatural transformation. By the parameter theorem for coends, $Z_T$ is an endofunctor of $\cc$ and $i_{X,Y}$ is
natural in $X$ and dinatural in~$Y$.

In \cite{BV3}, we construct an explicit a Hopf monad structure on $Z_T$, inherited from that of $T$. The Hopf monad $Z_T$ is
called the \emph{centralizer} of $T$.

Now, since $T$ preserves colimits (see~\cite[Remark 3.13]{BV2}) and so coends, $T(i)$ is a universal dinatural
transformation. Therefore we can define a natural transformation $\Omega\co TZ_T \to Z_TT$ by:
\begin{equation*}
\Omega_X T(i_{X,Y})= i_{T(X),T(Y)}\bigl(\ldual{\mu_Y}s^l_{T(Y)}T(\ldual{\mu_Y}) \otimes T_2(X,Y)\bigr)T_2(\ldual{T(Y)},X
\otimes Y),
\end{equation*}
where $\eta$ and $u$ the units of $T$ and $Z_T$ respectively, and $s^l$ is the left antipode of $T$.

\begin{thm}[\cite{BV3}]\label{thm-candistlaw}
$\Omega\co TZ_T \to Z_TT$ is a bijective comonoidal distributive law\footnote{A comonoidal distributive law between two Hopf
monads makes their composition a Hopf monad.}.
\end{thm}

The distributive law $\Omega$ is called the \emph{canonical distributive law} of $T$ over $Z_T$. Since~$\Omega$ is a
comonoidal distributive law, we get that $D_T=Z_T \circ_{\Omega} T$ is a Hopf monad on $\cc$ (whose underlying endofunctor
is $Z_T\circ T$). We call $D_T$ the \emph{double} of~$T$, as justified by the following theorem:

\begin{thm}[\cite{BV3}]\label{thm-doublable}
Let $T$ be a centralizable Hopf monad on an autonomous category~$\cc$. Then the forgetful functor $\uu\co \zz(T\ti\cc) \to
\cc$, given by $\bigl ((M,r),\sigma\bigr ) \mapsto M$, is monadic with monad the double $D_T$ of~$T$. Furthermore:
\begin{equation*}
R_{X,Y}=( u_{T(Y)} \otimes Z_T(\eta_X)\bigr)(\id_{T(Y)} \otimes i_{X,Y})(\lcoev_{Y} \otimes \id_X)
\end{equation*}
is a \Rt matrix for $D_T$, making the Hopf monad $D_T$ quasitriangular, and
\begin{equation*}
\zz(T\ti\cc) \cong D_T \ti \cc
\end{equation*}
as braided categories.
\end{thm}

\begin{rem}
Let $\cc$ be an autonomous category which is \emph{centralizable}, that is, such that the trivial Hopf monad $1_\cc$ is
centralizable. In that case, the centralizer $Z=Z_{1_\cc}$ and the double $D_{1_\cc}$ of $1_\cc$ coincide. Then, by
Theorem~\ref{thm-doublable}, $Z$ is a quasitriangular Hopf monad on $\cc$ such that $\zz(\cc) \cong Z \ti \cc$ as braided
category. In particular $\zz(\cc)$ is seen as the category of modules over a quasitriangular Hopf monad. In
Section~\ref{sect-centerfusion}, we explicitly describe $Z$ in terms of $\cc$ when $\cc$ is a fusion category.
\end{rem}

\begin{exa}\label{rem-doublable}
Let $H$ be a finite-dimensional Hopf algebra over a field $\kk$. Then the Hopf monad $T=?\otimes_\kk H$ on $\vect_\kk$ is
centralizable. We have: $Z_T=? \otimes_\kk H^*$ and so $D_T=?\otimes_\kk H \otimes_\kk H^*$. From
Theorem~\ref{thm-doublable}, the vector space $D(H)=H\otimes_\kk H^*$ inherits a quasitriangular Hopf algebra structure from
the quasitriangular Hopf monad $D_T$. In particular the algebra structure on $D(H)$ is a twist of that of $H \otimes H^*$ by
an isomorphism $H^* \otimes H \to H \otimes H^*$ coming from the distributive law $\Omega\co T Z_T \to Z_TT$. This
quasitriangular Hopf algebra $D(H)$ is precisely the \emph{Drinfeld double} of $H$. Furthermore, since
$T\ti\vect_\kk=\rep{H}$ and $D_T\ti\vect_\kk=\rep{D(H)}$, one recovers that $\zz(\rep{H}) \cong \rep{D(H)}$ as braided
categories.
\end{exa}

The previous example may be generalized to Hopf algebras in braided categories. Indeed, let $\cc$ be a braided category
which admits a coend:
\begin{equation*}
C=\int^{Y \in \cc} \ldual{Y} \otimes Y.
\end{equation*}
Recall $C$ is then a Hopf algebra in $\cc$ (see Section~\ref{sect-coend}). Let $A$ be a Hopf algebra in~$\cc$. Then the Hopf
monad $? \otimes A$ on $\cc$ is centralizable and we have:
\begin{equation*}
Z_{? \otimes A}=? \otimes \ldual{A} \otimes C, \qquad D_{? \otimes A}=? \otimes A \otimes \ldual{A} \otimes C.
\end{equation*}
From Theorem~\ref{thm-doublable}, we get that the object $D(A)=A \otimes \ldual{A} \otimes C$ is a quasitriangular Hopf
algebra in $\cc$, whose structure is inherited from the quasitriangular Hopf monad~$D_{? \otimes A}$. Here $D(A)$
quasitriangular means that there exists a \Rt matrix:
\begin{equation*}
R \co C \otimes C \to D(A) \otimes D(A)
\end{equation*}
verifying axioms generalizing the usual ones (when $\cc=\vect_\kk$, we have $C=\kk$ and $R\in D(H) \otimes D(H)$). This \Rt
matrix makes the category $\rep_\cc{D(A)}$ of right $D(A)$\ti modules (in $\cc$) braided, so that $\zz(\rep_\cc{A}) \cong
\rep_\cc{D(A)}$ as braided categories. We refer to~\cite{BV3} for more details.\\

\subsection{The coend of a category of modules over a Hopf monad}
Let $T$ be a centralizable Hopf monad on an autonomous category $\cc$. Denote $Z_T$ its centralizer, $i_{X,Y} \co
\ldual{T(Y)} \otimes X \otimes Y \to Z_T(X)$ its associated universal dinatural transformation, and $\Omega\co TZ_T \to
Z_TT$ the canonical distributive law of $T$ over $Z_T$. Then :
\begin{thm}[\cite{BV3}]\label{thm-coend}
The category $T\ti\cc$ of $T$-modules admits a coend, which is:
\begin{equation*}
\int^{(M,r) \in T\ti\cc}\!\!\!\!\!\!\!\!\!\!\!\!\!\!\!\!\!\!\!\!\! \ldual{(M,r)} \otimes
(M,r)=\bigl(Z_T(\un),Z_T(T_0)\Omega_\un\bigr),
\end{equation*}
with $I_{(M,r)}=i_{\un,M}(\ldual{r} \otimes \id_M)\co \ldual{(M,r)} \otimes (M,r) \to
\bigl(Z_T(\un),Z_T(T_0)\Omega_\un\bigr)$ as universal dinatural transformation.
\end{thm}

Note that if $T$ is furthermore quasitriangular, then $T\ti\cc$ is braided and so the coend
$\bigl(Z_T(\un),Z_T(T_0)\Omega_\un\bigr)$ is a Hopf algebra in $T\ti\cc$.

\begin{rem}
If we apply this to the double $D_T$ of a centralizable Hopf monad, which we suppose to be itself centralizable, we get an
explicit description of the coend of the braided category $D_T \ti \cc \cong \zz(T\ti\cc)$ in terms of the monad $T$. Using
this description and Hopf diagrams, we hence have a way of computing the quantum invariants of 3-manifolds defined using
$\zz(T\ti\cc)$, see Section~\ref{sect-Hdiag}. In the next section, we explicit the case where $\cc$ is a spherical fusion
category and $T=1_\cc$.
\end{rem}

\section{Reshetikhin-Turaev invariants from categorical centers}\label{sect-RT}

In this section, we treat in details the case of the center $\zz(\cc)$ of a spherical fusion category $\cc$. This leads to
an explicit algorithm for computing Reshetikhin-Turaev-like invariants defined using $\zz(\cc)$ in terms of $\cc$.

\subsection{On the center of a fusion category}\label{sect-centerfusion} Fix a commutative ring $\kk$. Let $\cc$ be a fusion
category over $\kk$ (see Section~\ref{sect-fusion}). Then the trivial Hopf monad $1_\cc$ is centralizable. Its centralizer $Z=Z_{1_\cc}$ is:
\begin{equation*}
Z(X)=\bigoplus_{i \in I} \ldual{V}_i \otimes X \otimes V_i,
\end{equation*}
with associated dinatural transformation $i_{X,Y} \co \ldual{Y} \otimes X \otimes Y \to Z(X)$ given by:
\begin{equation*}
i_{X,Y}=\sum_{\substack{i\in I \\ 1 \leq \alpha \leq N_Y^i}} \ldual{q}_Y^{i,\alpha}\otimes \id_X \otimes p_Y^{i,\alpha}.
\end{equation*}

The double of $1_\cc$ is $D_{1_\cc}=Z \circ 1_\cc=Z$. Hence $Z$ is a quasitriangular Hopf monad and $\zz(\cc) \cong Z\ti\cc$
as braided categories. Furthermore, if $\cc$ is spherical, then $Z$ is a ribbon Hopf monad (and so $\zz(\cc)$ is ribbon).

\begin{figure}[t]
   \begin{center}
          $\displaystyle Z_2(X,Y)=\sum_{i \in I}$\,
 \psfrag{o}[Bc][Bc]{\scalebox{.55}{$V_i$}}
 \psfrag{n}[Bc][Bc]{\scalebox{.55}{$X$}}
 \psfrag{u}[Bc][Bc]{\scalebox{.55}{$Y$}}
 \psfrag{c}[Bc][Bc]{\scalebox{.55}{$\ldual{V}_i$}}
\scaleraisedraw{.4}{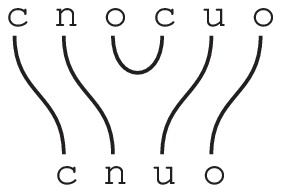} , \qquad $\displaystyle Z_0=\sum_{i \in I}$\,
 \psfrag{o}[Bc][Bc]{\scalebox{.55}{$V_i$}}
 \psfrag{c}[Bc][Bc]{\scalebox{.55}{$\ldual{V}_i$}}
\scaleraisedraw{.6}{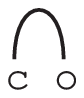} , \\[.5em]
$\displaystyle \mu_X=\hspace*{-.2cm}\sum_{\substack{i,j,k \in I \\ 1 \leq \alpha \leq N_{i,j}^k }}$\,
 \psfrag{u}[Bc][Bc]{\scalebox{.55}{$V_i$}}
 \psfrag{r}[Bc][Bc]{\scalebox{.55}{$V_j$}}
 \psfrag{a}[Bc][Bc]{\scalebox{.55}{$V_k$}}
 \psfrag{o}[Bc][Bc]{\scalebox{.55}{$\ldual{V}_k$}}
 \psfrag{c}[Bc][Bc]{\scalebox{.55}{$\ldual{V}_j$}}
 \psfrag{s}[Bc][Bc]{\scalebox{.55}{$\ldual{V}_i$}}
 \psfrag{n}[Bc][Bc]{\scalebox{.55}{$X$}}
 \psfrag{p}[c][c]{\scalebox{.7}{$\ldual{q}^{k,\alpha}_{i,j}$}}
 \psfrag{q}[c][c]{\scalebox{.7}{$p^{k,\alpha}_{i,j}$}}
\scaleraisedraw{.5}{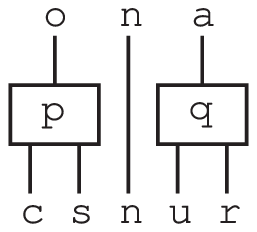} , \qquad $\displaystyle \eta_X=$
 \psfrag{n}[Bc][Bc]{\scalebox{.55}{$V_0$}}
 \psfrag{c}[Bc][Bc]{\scalebox{.55}{$X$}}
 \psfrag{o}[Bc][Bc]{\scalebox{.55}{$\ldual{V}_0$}}
\scaleraisedraw{.35}{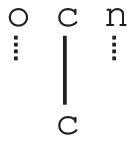} , \\[.5em]
$\displaystyle s^l_X=\sum_{i\in I}$
 \psfrag{u}[Bc][Bc]{\scalebox{.55}{$\lldual{V}_i$}}
 \psfrag{a}[Bc][Bc]{\scalebox{.55}{$V_{\ldual{i}}$}}
 \psfrag{c}[Bc][Bc]{\scalebox{.55}{$\ldual{V}_{\ldual{i}}$}}
 \psfrag{s}[Bc][Bc]{\scalebox{.55}{$\ldual{V}_i$}}
 \psfrag{n}[Bc][Bc]{\scalebox{.55}{$\ldual{X}$}}\,
\scaleraisedraw{.35}{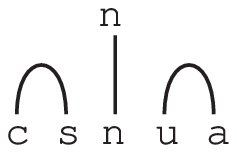} , \qquad $\displaystyle s^r_X=\sum_{i\in I}$
 \psfrag{u}[Bc][Bc]{\scalebox{.55}{$V_i$}}
 \psfrag{a}[Bc][Bc]{\scalebox{.55}{$V_{\rdual{i}}$}}
 \psfrag{c}[Bc][Bc]{\scalebox{.55}{$\ldual{V}_{\rdual{i}}$}}
 \psfrag{s}[Bc][Bc]{\scalebox{.55}{$\rdual{V}_i$}}
 \psfrag{n}[Bc][Bc]{\scalebox{.55}{$\rdual{X}$}}\,
\scaleraisedraw{.35}{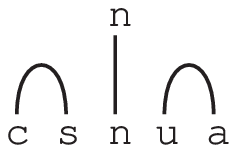} , \\[.5em]
$\displaystyle R_{X,Y}=\hspace*{-.2cm}\sum_{\substack{i \in I \\ 1 \leq \alpha \leq N_Y^i }}$\hspace*{-.2cm}
 \psfrag{u}[Bc][Bc]{\scalebox{.55}{$V_i$}}
 \psfrag{a}[Bc][Bc]{\scalebox{.55}{$V_0$}}
 \psfrag{c}[Bc][Bc]{\scalebox{.55}{$\ldual{V}_0$}}
 \psfrag{s}[Bc][Bc]{\scalebox{.55}{$\ldual{V}_i$}}
 \psfrag{n}[Bc][Bc]{\scalebox{.55}{$X$}}
 \psfrag{o}[Bc][Bc]{\scalebox{.55}{$Y$}}
 \psfrag{p}[c][c]{\scalebox{.7}{$p_Y^i$}}
 \psfrag{q}[c][c]{\scalebox{.7}{$q_Y^i$}}
\scaleraisedraw{.5}{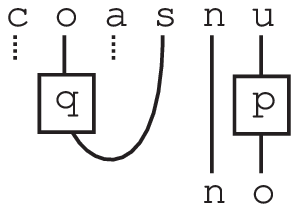}  , \qquad $\displaystyle \theta_X=\hspace*{-.2cm}\sum_{\substack{i \in I \\ 1 \leq \alpha \leq
N_X^i }}$
 \psfrag{u}[Bc][Bc]{\scalebox{.55}{$V_{\ldual{i}}$}}
 \psfrag{s}[Bc][Bc]{\scalebox{.55}{$\ldual{V}_{\ldual{i}}$}}
 \psfrag{o}[Bc][Bc]{\scalebox{.55}{$X$}}
 \psfrag{p}[c][c]{\scalebox{.7}{$p_X^i$}}
 \psfrag{q}[c][c]{\scalebox{.7}{$q_X^i$}}
 \psfrag{g}[c][c]{\scalebox{.7}{$\phi_{V_i}$}}
\scaleraisedraw{.5}{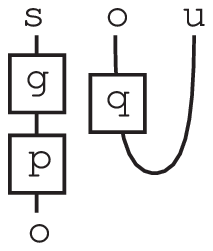}  .
   \end{center}
     \caption{Structural morphisms of $Z$}
     \label{morphZ}
\end{figure}

The structural morphisms of $Z$ can be described only in terms of the category~$\cc$, that is, only using the $p,q$'s (see Section~\ref{sect-fusion}), the duality morphisms, and the sovereign structure $\phi_X \co X \to
\lldual{X}$. They are depicted in Figure~\ref{morphZ}. The dotted lines in the figure represent $\id_{V_0}=\id_\un$ and can
be removed without changing the morphisms. We depicted them in order to remember which factor of~$Z(X)$ is concerned. To
simplify the reading, we denote $A_{V_{i_1} \otimes \cdots \otimes V_{i_n}}$ by $A_{i_1, \dots, i_n}$ for $A=p^{i,\alpha}$,
$q^{i,\alpha}$, or~$N^i$.

Using the Maschke theorem for Hopf monad's which characterize semisimplicity (see \cite[Theorem 6.5]{BV2}), we have:
\begin{prop}\cite{BV4}\label{propZss}
Let $\cc$ be a spherical fusion category. Then the (ribbon) Hopf monad $Z$ is semisimple if and only if $\dim \cc$ is
invertible.
\end{prop}

Since $\zz(\cc) \cong Z\ti\cc$, a direct consequence of Proposition~\ref{propZss} is then:
\begin{cor}\label{cor-centerfusion}
Let $\cc$ be a spherical fusion category over an algebraic closed field. Assume $\dim \cc$ is invertible. Then $\zz(\cc)$ is
a ribbon fusion category.
\end{cor}

\subsection{The coend of the center of a fusion category}\label{sect-fusioncoend}
Let us describe the structure of the coend $(C,r)$ of $Z\ti \cc\cong \zz(\cc)$, where $\cc$ is a spherical fusion category.
Recall that $C$ is an object of $\cc$ and $r\co Z(C) \to C$ is an action of $Z$ on $C$. From Theorem~\ref{thm-coend}, we
get:
\begin{equation*}
C=\bigoplus_{j \in I} \ldual{Z(V_j)} \otimes V_j= \bigoplus_{i,j \in I} \ldual{V}_i \otimes \ldual{V}_j \otimes \lldual{V}_i
\otimes V_j.
\end{equation*}
Note that an immediate consequence of this is: $\dim \zz(\cc) =(\dim \cc)^2$.

\begin{figure}[t]
   \begin{center}
       $\displaystyle \Delta_C=\hspace*{-.2cm}\sum_{\substack{i,j,k,m,n \in I \\ 1 \leq \alpha \leq N_{k,m}^k \\ 1 \leq
\beta \leq N_{\ldual{k},j,k}^n}}$
 \psfrag{u}[Bc][Bc]{\scalebox{.55}{$\ldual{V}_i$}}
 \psfrag{r}[Bc][Bc]{\scalebox{.55}{$\ldual{V}_j$}}
 \psfrag{i}[Bc][Bc]{\scalebox{.55}{$\lldual{V}_i$}}
 \psfrag{b}[Bc][Bc]{\scalebox{.55}{$V_j$}}
 \psfrag{o}[Bc][Bc]{\scalebox{.55}{$\ldual{V}_m$}}
 \psfrag{c}[Bc][Bc]{\scalebox{.55}{$\ldual{V}_n$}}
 \psfrag{s}[Bc][Bc]{\scalebox{.55}{$\lldual{V}_m$}}
 \psfrag{n}[Bc][Bc]{\scalebox{.55}{$V_n$}}
 \psfrag{w}[Bc][Bc]{\scalebox{.55}{$\ldual{V}_k$}}
 \psfrag{e}[Bc][Bc]{\scalebox{.55}{$\ldual{V}_j$}}
 \psfrag{z}[Bc][Bc]{\scalebox{.55}{$\lldual{V}_k$}}
 \psfrag{p}[c][c]{\scalebox{.7}{$\ldual{p}^{n,\beta}_{\ldual{k},j,k}$}}
 \psfrag{q}[c][c]{\scalebox{.7}{$\ldual{q}^{n,\beta}_{\ldual{k},j,k}$}}
 \psfrag{a}[c][c]{\scalebox{.7}{$\ldual{p}^{i,\alpha}_{k,m}$}}
 \psfrag{x}[c][c]{\scalebox{.7}{$\lldual{q}^{i,\alpha}_{k,m}$}}
\scaleraisedraw{.6}{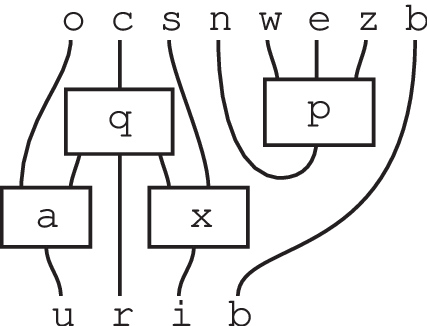} , \qquad $\displaystyle\varepsilon_C= \sum_{j \in I}$ \,
 \psfrag{o}[Bc][Bc]{\scalebox{.55}{$\ldual{V}_0$}}
 \psfrag{c}[Bc][Bc]{\scalebox{.55}{$\ldual{V}_j$}}
 \psfrag{n}[Bc][Bc]{\scalebox{.55}{$\lldual{V}_0$}}
 \psfrag{e}[Bc][Bc]{\scalebox{.55}{$V_j$}}
\scaleraisedraw{.4}{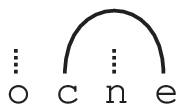} ,\\[.5em]
$\displaystyle S_C=\hspace*{-.6cm}\sum_{\substack{i,j,k,l \in I \\ 1 \leq \alpha \leq N^{\ldual{i}}_{j,k,\rdual{j}} \\ 1
\leq \beta \leq N^l_{\ldual{j},\ldual{i}, \ldual{j}, \lldual{i}, j}}}$
 \psfrag{x}[c][c]{\scalebox{.7}{$p^{l,\beta}_{\ldual{j},\ldual{i}, \ldual{j}, \lldual{i}, j}$}}
 \psfrag{a}[c][c]{\scalebox{.7}{$\ldual{q}^{l,\beta}_{\ldual{j},\ldual{i}, \ldual{j}, \lldual{i}, j}$}}
 \psfrag{p}[c][c]{\scalebox{.7}{$\ldual{p}^{\ldual{i},\alpha}_{j,k,\rdual{j}}$}}
 \psfrag{q}[c][c]{\scalebox{.7}{$\lldual{q}^{\ldual{i},\alpha}_{j,k,\rdual{j}}$}}
 \psfrag{o}[Bc][Bc]{\scalebox{.55}{$\ldual{V}_k$}}
 \psfrag{c}[Bc][Bc]{\scalebox{.55}{$\ldual{V}_l$}}
 \psfrag{s}[Bc][Bc]{\scalebox{.55}{$\lldual{V}_k$}}
 \psfrag{n}[Bc][Bc]{\scalebox{.55}{$V_l$}}
\scaleraisedraw{.5}{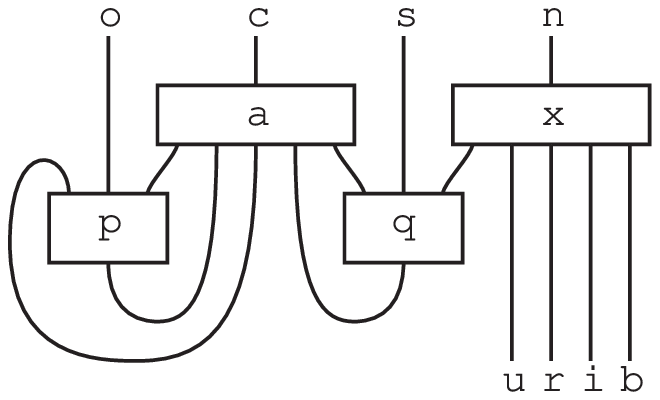} ,\\[.5em]
$\displaystyle \omega_C=\hspace*{-.2cm}\sum_{\substack{i,j,k,l \in I \\ 1 \leq \alpha \leq N^{\ldual{k}}_{\ldual{i},j,i} \\
1 \leq \beta \leq N^i_{\ldual{k},\ldual{l},\lldual{k}}}}$
 \psfrag{x}[c][c]{\scalebox{.7}{$p^{\ldual{k},\alpha}_{\ldual{i},j,i}$}}
 \psfrag{a}[c][c]{\scalebox{.7}{$q^{i,\beta}_{\ldual{k},\ldual{l},\lldual{k}}$}}
 \psfrag{p}[c][c]{\scalebox{.7}{$\ldual{q}^{\ldual{k},\alpha}_{\ldual{i},j,i}$}}
 \psfrag{q}[c][c]{\scalebox{.7}{$p^{i,\beta}_{\ldual{k},\ldual{l},\lldual{k}}$}}
 \psfrag{o}[Bc][Bc]{\scalebox{.55}{$\ldual{V}_i$}}
 \psfrag{c}[Bc][Bc]{\scalebox{.55}{$\ldual{V}_j$}}
 \psfrag{s}[Bc][Bc]{\scalebox{.55}{$\lldual{V}_i$}}
 \psfrag{n}[Bc][Bc]{\scalebox{.55}{$V_j$}}
 \psfrag{u}[Bc][Bc]{\scalebox{.55}{$\ldual{V}_k$}}
 \psfrag{r}[Bc][Bc]{\scalebox{.55}{$\ldual{V}_l$}}
 \psfrag{i}[Bc][Bc]{\scalebox{.55}{$\lldual{V}_k$}}
 \psfrag{b}[Bc][Bc]{\scalebox{.55}{$V_l$}} \;
\scaleraisedraw{.55}{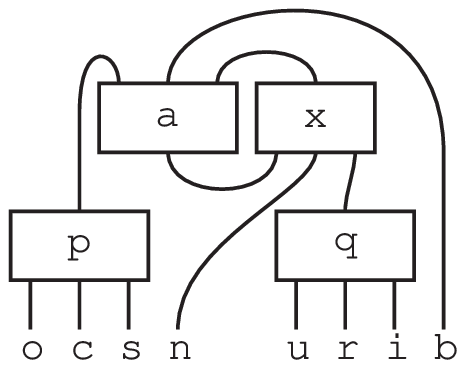} ,\\[.5em]
$\displaystyle\theta_C^+=\sum_{i \in I}$\;
 \psfrag{a}[c][c]{\scalebox{.7}{$\phi_{\ldual{V}_i}$}}
 \psfrag{o}[Bc][Bc]{\scalebox{.55}{$\ldual{V}_i$}}
 \psfrag{c}[Bc][Bc]{\scalebox{.55}{$\ldual{V}_i$}}
 \psfrag{n}[Bc][Bc]{\scalebox{.55}{$\lldual{V}_i$}}
 \psfrag{e}[Bc][Bc]{\scalebox{.55}{$V_i$}}
\scaleraisedraw{.47}{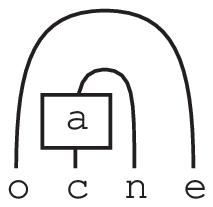}\, , \qquad $\displaystyle\theta_C^-= \sum_{i \in I}$\;
 \psfrag{a}[c][c]{\scalebox{.7}{$\phi_{\ldual{V}_i}$}}
 \psfrag{o}[Bc][Bc]{\scalebox{.55}{$\ldual{V}_i$}}
 \psfrag{c}[Bc][Bc]{\scalebox{.55}{$\ldual{V}_{\ldual{i}}$}}
 \psfrag{n}[Bc][Bc]{\scalebox{.55}{$\lldual{V}_i$}}
 \psfrag{e}[Bc][Bc]{\scalebox{.55}{$V_{\ldual{i}}$}}
\scaleraisedraw{.5}{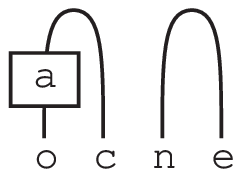} , \\[.5em]
$\displaystyle\tau_{C,C}=\hspace*{-.1cm}\sum_{\substack{i,j,k,l,a,b,m \in I \\
1 \leq \alpha \leq N^m_{\ldual{k},\ldual{l},\lldual{k},l} \\
\\ 1 \leq \beta \leq N^b_{\ldual{m},j,m}\\ 1 \leq \gamma \leq N^i_{m,a,\rdual{m}}}}$ \hspace*{-.4cm}
 \psfrag{x}[c][c]{\scalebox{.7}{$\ldual{p}^{i,\gamma}_{m,a,\rdual{m}}$}}
 \psfrag{v}[c][c]{\scalebox{.7}{$\lldual{q}^{i,\gamma}_{m,a,\rdual{m}}$}}
 \psfrag{p}[c][c]{\scalebox{.7}{$p^{m,\alpha}_{\ldual{k},\ldual{l},\lldual{k},l}$}}
 \psfrag{q}[c][c]{\scalebox{.7}{$q^{m,\alpha}_{\ldual{k},\ldual{l},\lldual{k},l}$}}
 \psfrag{o}[c][c]{\scalebox{.7}{$\ldual{q}^{b,\beta}_{\ldual{m},j,m}$}}
 \psfrag{a}[c][c]{\scalebox{.7}{$p^{b,\beta}_{\ldual{m},j,m}$}}
 \psfrag{w}[Bc][Bc]{\scalebox{.55}{$\ldual{V}_i$}}
 \psfrag{c}[Bc][Bc]{\scalebox{.55}{$\ldual{V}_j$}}
 \psfrag{s}[Bc][Bc]{\scalebox{.55}{$\lldual{V}_i$}}
 \psfrag{n}[Bc][Bc]{\scalebox{.55}{$V_j$}}
 \psfrag{u}[Bc][Bc]{\scalebox{.55}{$\ldual{V}_k$}}
 \psfrag{r}[Bc][Bc]{\scalebox{.55}{$\ldual{V}_l$}}
 \psfrag{i}[Bc][Bc]{\scalebox{.55}{$\lldual{V}_k$}}
 \psfrag{b}[Bc][Bc]{\scalebox{.55}{$V_l$}} \;
 \psfrag{t}[Bc][Bc]{\scalebox{.55}{$\ldual{V}_a$}}
 \psfrag{l}[Bc][Bc]{\scalebox{.55}{$\ldual{V}_b$}}
 \psfrag{h}[Bc][Bc]{\scalebox{.55}{$\lldual{V}_a$}}
 \psfrag{f}[Bc][Bc]{\scalebox{.55}{$V_b$}} \;
\scaleraisedraw{.6}{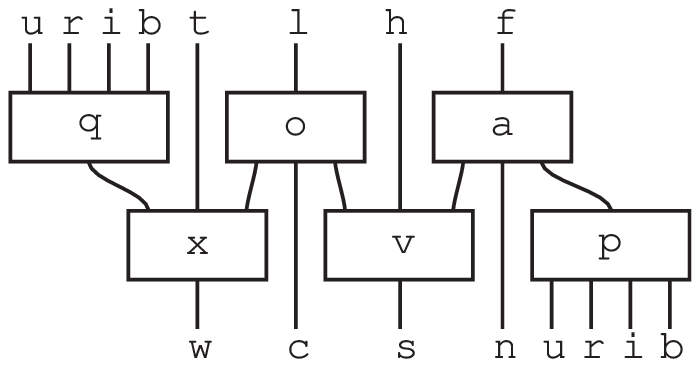} , \\[.5em]
$\displaystyle\Lambda=\sum_{j \in I} \dim (V_j) $\;
 \psfrag{a}[c][c]{\scalebox{.7}{$\phi^{-1}_{V_j}$}}
 \psfrag{o}[Bc][Bc]{\scalebox{.55}{$\ldual{V}_0$}}
 \psfrag{c}[Bc][Bc]{\scalebox{.55}{$\ldual{V}_j$}}
 \psfrag{n}[Bc][Bc]{\scalebox{.55}{$\lldual{V}_0$}}
 \psfrag{e}[Bc][Bc]{\scalebox{.55}{$V_j$}}
\scaleraisedraw{.45}{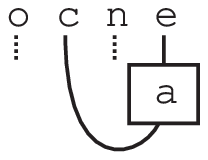} .
   \end{center}
     \caption{Structural morphisms of the coend of $\zz(\cc)$}
     \label{morphcoendcenter}
\end{figure}

The structural morphisms of $C$ can be expressed using only the category $\cc$. Those needed to represent Hopf diagrams are
depicted in Figure~\ref{morphcoendcenter}.

\begin{thm}[\cite{BV4}]\label{thm-coint-coend}
The morphism $\Lambda\co \un \to C$ of Figure~\ref{morphcoendcenter} is a $S_C$-invariant integral of the coend of
$Z\ti\cc\cong \zz(\cc)$.
\end{thm}

Following~\cite{Lyu2}, a braided category $\bb$ is said to be \emph{modular} if it admits a coend~$C$ whose Hopf pairing
$\omega_C\co C \otimes C \to \un$ is non-degenerate (meaning there exists $\sigma\co \un \to C \otimes C$ such that
$(\omega_C \otimes \id_C)(\id_C \otimes \sigma)=\id_C=(\id_C \otimes \omega_C)(\sigma \otimes \id_C)$). Note this extends
the usual notion of modularity to the non-semisimple case (when~$\bb$ is a ribbon fusion category, $\bb$ is modular in the
above sense if and only if the $S$\ti matrix is invertible).
\begin{cor}[\cite{BV4}]
The center of a spherical fusion category is modular.
\end{cor}

\begin{rem}\label{rem-algclosedmod}
Let $\cc$ a spherical fusion category over an algebraic closed field such that $\dim \cc$ is invertible. Then by
Corollary~\ref{cor-centerfusion} and Theorem~\ref{thm-coint-coend}, we get that the center $\zz(\cc)$ of $\cc$ is a modular
ribbon fusion category. This last result was first shown in~\cite{Mueg} using different method.
\end{rem}

\subsection{Computing $\mathrm{RT}_{\zz(\cc)}(M^3)$ from $\cc$}\label{soussect-RT}
Let $\cc$ be a spherical fusion category over a commutative ring $\kk$. As explained in Sections~\ref{sect-centerfusion}
and~\ref{sect-fusioncoend}, the center $Z(\cc)$ of $\cc$ is a ribbon category which admits a coend $C$.

The integral~$\Lambda$ of $C$ is then a normalizable Kirby element since $\theta_C^+ \Lambda=1_\kk$ and
$\theta_C^-\Lambda=1_\kk$. Hence the invariant $\tau_{\zz(\cc)}(M,\Lambda)$ of $3$\ti manifolds (see
Section~\ref{sect-kirby}).

Furthermore, since we have an explicit description of the structural morphisms of the coend $C$ (see
Figures~\ref{morphcoendcenter}), we have a way to compute this invariant by using Hopf diagrams (see
Section~\ref{sect-Hdiag}). For example, we have:
\begin{equation*}
\tau_{\zz(\cc)}(S^3;\Lambda)=1 \quad \text{and} \quad \tau_{\zz(\cc)}(S^2 \times S^1;\Lambda)=\dim \cc.
\end{equation*}

Note that the invariant $\tau_{\zz(\cc)}(M,\Lambda)$ is well-defined even if $\dim \cc$ is not invertible. When $\dim \cc$
is invertible and $\kk$ is an algebraic closed field (so that $\zz(\cc)$ is a modular fusion category, see
Remark~\ref{rem-algclosedmod}), the invariant $\tau_{\zz(\cc)}(M,\Lambda)$ equals to the Reshetikhin-Turaev invariant
$\mathrm{RT}_{\zz(\cc)}(M)$ (up to a different normalization, see Remark~\ref{rem-invRT}). Hence we get a way to compute
$\mathrm{RT}_{\zz(\cc)}(M)$ in terms of the structural morphisms of $\cc$ (recall one
cannot use the original algorithm of Reshetikhin-Turaev since the simple objects of $\zz(\cc)$ are unknown in general).

\end{document}